\documentclass[12pt,a4paper]{article}%
\usepackage[utf8]{inputenc}
\usepackage{amsmath}
\usepackage{amsfonts}
\usepackage{amssymb}
\usepackage{graphicx}
\usepackage[space]{grffile}
\usepackage{hyperref}
\usepackage{layout}
\usepackage{longtable}
\usepackage{xcolor}%
\providecommand{\U}[1]{\protect\rule{.1in}{.1in}}
\newtheorem{theorem}{Theorem}
\newtheorem{theorem*}{Example}

\newtheorem{conjecture}[theorem]{Conjecture}
\newtheorem{corollary}[theorem]{Corollary}

\newtheorem{definition}[theorem]{Definition}

\newtheorem{lemma}[theorem]{Lemma}

\newtheorem{proposition}[theorem]{Proposition}
\newtheorem{remark}[theorem]{Observation}

\newenvironment{proof}[1][Proof]{\noindent\textbf{#1.} }{\ \hfill \rule{0.5em}{0.5em}\bigskip}
\graphicspath{{D:/Dropbox/Riste-Sedlar/Edge colorings/Paper2/Slike/}}
\begin{document}

\title{Normal $5$-edge-coloring of some snarks superpositioned by Flower snarks}
\author{Jelena Sedlar$^{1,3}$,\\Riste \v Skrekovski$^{2,3}$ \\[0.3cm] {\small $^{1}$ \textit{University of Split, FGAG, 21000 Split,
Croatia}}\\[0.1cm] {\small $^{2}$ \textit{University of Ljubljana, FMF, 1000 Ljubljana,
Slovenia }}\\[0.1cm] {\small $^{3}$ \textit{Faculty of Information Studies, 8000 Novo
Mesto, Slovenia }}\\[0.1cm] }
\maketitle

\begin{abstract}
An edge $e$ is normal in a proper edge-coloring of a cubic graph $G$ if the
number of distinct colors on four edges incident to $e$ is $2$ or $4.$ A
normal edge-coloring of $G$ is a proper edge-coloring in which every edge of
$G$ is normal. The Petersen Coloring Conjecture is equivalent to stating that
every bridgeless cubic graph has a normal $5$-edge-coloring. Since every
$3$-edge-coloring of a cubic graph is trivially normal, it is sufficient to
consider only snarks to establish the conjecture. In this paper, we consider a
class of superpositioned snarks obtained by choosing a cycle $C$ in a snark
$G$ and superpositioning vertices of $C$ by one of two simple supervertices
and edges of $C$ by superedges $H_{x,y}$, where $H$ is any snark and $x,y$ any
pair of nonadjacent vertices of $H.$ For such superpositioned snarks, two
sufficient conditions are given for the existence of a normal $5$%
-edge-coloring. The first condition yields a normal $5$-edge-coloring for all
hypohamiltonian snarks used as superedges, but only for some of the possible
ways of connecting them. In particular, since the Flower snarks are
hypohamiltonian, this consequently yields a normal $5$-edge-coloring for many
snarks superpositioned by the Flower snarks. The second sufficient condition
is more demanding, but its application yields a normal $5$-edge-colorings for
all superpositions by the Flower snarks. The same class of snarks is
considered in \textit{[S. Liu, R.-X. Hao, C.-Q. Zhang, Berge--Fulkerson
coloring for some families of superposition snarks, Eur. J. Comb. 96 (2021)
103344]} for the Berge-Fulkerson conjecture. Since we established that this
class has a Petersen coloring, this immediately yields the result of the above
mentioned paper.

\end{abstract}

\textit{Keywords:} normal edge-coloring; cubic graph; snark; superposition;
Petersen Coloring Conjecture.

\textit{AMS Subject Classification numbers:} 05C15

\section{Introduction}

In this paper we consider simple cubic bridgeless graphs. Let $G$ be a graph
with the set of vertices $V(G)$ and the set of edges $E(G).$ For a vertex
$v\in V(G)$ we denote by $\partial_{G}(v)$ the set of edges in $G$ which are
incident to $v.$ Let $G$ and $H$ be a pair of cubic graphs. A $H$%
\emph{-coloring} of $G$ is any mapping $\phi:E(G)\rightarrow E(H)$ such that
for every $v\in V(G)$ there is $w\in V(H)$ such that $\phi(\partial
_{G}(v))=\partial_{H}(w).$ If $G$ admits an $H$-coloring, then we write
$H\prec G.$

Let us denote by $P_{10}$ the Petersen graph. If $H=P_{10}$ and $H\prec G,$
then we say that $G$ has the Petersen coloring. The following conjecture was
stated by Jaeger \cite{Jaeger1988}.

\begin{conjecture}
[Petersen Coloring Conjecture]\label{Con_petersen}Let $G$ be a bridgeless
cubic graph. Then $P_{10}\prec G.$
\end{conjecture}

It is known that the Petersen Coloring Conjecture implies the classical
Berge-Fulkerson Conjecture and the (5,2)-cycle-cover conjecture. Thus, it
seems difficult to prove. It motivated several variants and different
approaches to it \cite{Mkrtchjan2013, Riste2020, Samal2011}, one such approach
are normal colorings.

\paragraph{Normal colorings.}

A \emph{(proper) }$k$\emph{-edge-coloring} of a graph $G$ is any mapping
$\sigma:E(G)\rightarrow\{1,\ldots,k\}$ such that any pair of adjacent edges of
$G$ receives distinct colors by $\sigma$. Let $\sigma$ be a $k$-edge-coloring
of $G$ and $v\in V(G).$ Then we denote $\sigma(v)=\{\sigma(e):e\in\partial
_{G}(v)\}.$

\begin{definition}
Let $G$ be a bridgeless cubic graph, $\sigma$ a proper edge-coloring of $G$
and $uv$ an edge of $G.$ We say the edge $uv$ is \emph{poor} (resp.
\emph{rich}) if $\left\vert \sigma(u)\cup\sigma(v)\right\vert =3$ (resp.
$\left\vert \sigma(u)\cup\sigma(v)\right\vert =5$).
\end{definition}

A proper edge-coloring of a cubic graph $G$ is \emph{normal}, if every edge of
$G$ is poor or rich. Normal edge-colorings were introduced by Jaeger in
\cite{Jaeger1985}. Notice that a normal coloring in which every edge is poor
must be a $3$-edge-coloring. On the other hand, a normal coloring in which
every edge is rich is called a \emph{strong edge-coloring}. The \emph{normal
chromatic index} of a cubic graph $G,$ denoted by $\chi_{N}^{\prime}(G),$ is
the smallest number $k$ such that $G$ admits a normal $k$-edge-coloring.
Notice that $\chi_{N}^{\prime}(G)$ is at least $3$ and it never equals $4$.
The following result was established in \cite{Jaeger1985}.

\begin{proposition}
\label{Prop_equivalentNormalPetersen}For a cubic graph $G$, it holds that
$P_{10}\prec G$ if and only if $G$ admits a normal $5$-edge-coloring.
\end{proposition}

Proposition \ref{Prop_equivalentNormalPetersen} implies that Conjecture
\ref{Con_petersen} can be restated as follows.

\begin{conjecture}
\label{Con_normal}Let $G$ be a bridgeless cubic graph. Then $\chi_{N}^{\prime
}(G)\leq5.$
\end{conjecture}

\noindent Since any proper $3$-edge-coloring of a cubic graph $G$ is a normal
edge-coloring in which every edge is poor, Conjecture \ref{Con_normal}
obviously holds for every $3$-edge colorable graph. The smallest known upper
bound for $\chi_{N}^{\prime}(G)$ of a bridgeless cubic graph is $7$
\cite{BilkovaHana, Mazzucuolo2020normal}, and for several graph classes the
bound is lowered to $6$ \cite{Mazzucuolo2020normal6}. According to Vizing's
theorem, every cubic graph is either $3$-edge colorable or $4$-edge colorable,
so to prove Conjecture \ref{Con_normal} it remains to establish that it holds
for all bridgeless cubic graphs which are not $3$-edge colorable.

\paragraph{Superpositioning snarks.}

Cubic graphs which are not $3$-edge colorable are usually considered under the
name of snarks \cite{MazzuoccoloStephen,NedelaSkovieraSurvey}. In order to
avoid trivial cases, the definition of snark usually includes some additional
requirements on the connectivity and the girth. Since in this paper such
requirements are not essential, we will go with the broad definition of a
\emph{snark} being any bridgeless cubic graph which is not $3$-edge colorable.
The existence of a normal $5$-edge-coloring for some families of snarks has
already been established, see for example \cite{MazzuccoloLupekhine,
Hagglund2014}.

A general method for obtaining snarks is the superposition \cite{Adelson1973,
Descartes1946, KocholSuperposition, Kochol2, MacajovaRevisited}. We will focus
our attention to a class of snarks obtained by some particular superpositions
of edges and vertices of a cycle in a snark. So, let us first introduce the
method of superposition.

\begin{definition}
A \emph{multipole} $M=(V,E,S)$ is a triple which consists of a set of vertices
$V=V(M)$, a set of edges $E=E(M)$, and a set of semiedges $S=S(M)$. A semiedge
is incident either to one vertex or to another semiedge in which case a pair
of incident semiedges forms a so called \emph{isolated edge} within the multipole.
\end{definition}

\begin{figure}[h]
\begin{center}%
\begin{tabular}
[c]{ll}%
\begin{tabular}
[c]{c}%
\includegraphics[scale=0.7]{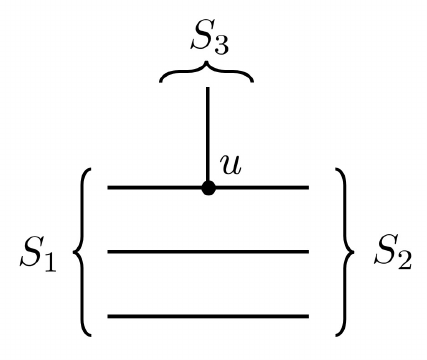}\\
$A$%
\end{tabular}
&
\begin{tabular}
[c]{c}%
\includegraphics[scale=0.7]{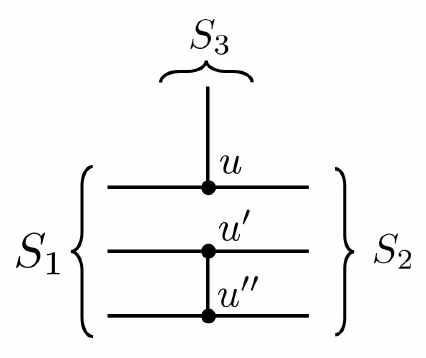}\\
$A^{\prime}$%
\end{tabular}
\end{tabular}
\end{center}
\caption{Multipoles $A$ and $A^{\prime}$. Notice that $A$ has two isolated
edges. Both $A$ and $A^{\prime}$ are $7$-poles, and if their set of semiedges
is divided into three connectors $S_{1},$ $S_{2}$ and $S_{3}$ as shown in the
figure, then $A$ and $A^{\prime}$ are supervertices.}%
\label{Fig_superVert}%
\end{figure}

For an illustration of a multipole see Figure \ref{Fig_superVert}. All
multipoles we consider here are \emph{cubic}, i.e. such that every vertex is
incident to precisely three edges or semiedges. Notice that the notion of
normal coloring can be naturally extended to multipoles. Namely, a
\emph{(proper) }$k$\emph{-edge-coloring} of a multipole $M=(V,E,S)$ is any
mapping $\sigma:E(M)\cup E(S)\rightarrow\{1,\ldots,k\}$ such that colors of
edges and semiedges which are incident to a same vertex are pairwise distinct.
Further, for a cubic multipole $M$ we define a \emph{normal coloring} as a
proper edge-coloring such that every edge is rich or poor. Notice that the
definition of a normal coloring of a multipole does not have any requirements
on the colors of semiedges. Finally, we need the notion of a Kempe chain in a
cubic multipole. So, let $\sigma$ be a normal coloring of a cubic multipole
$M=(V,E,S)$ and let $P=s_{1}e_{2}\cdots e_{p-1}s_{p}$, for $p>1,$ $e_{i}\in E$
and $s_{i}\in S,$ be a path in $M$ starting and ending with a semiedge. We say
that $P$ is a \emph{Kempe }$(i,j)$\emph{-chain} with respect to $\sigma,$ if
(semi)edges of $P$ are alternatively colored by $i$ and $j$.

Next, we wish to define special kinds of multipoles which will be of interest
to us. We start with the definition of $k$\emph{-pole} which is defined as any
multipole $M$ with $\left\vert S(M)\right\vert =k$. The set of semiedges of a
$k$-pole $M$ can be divided into several so called connectors, i.e. we define
a $(k_{1},\ldots,k_{n})$\emph{-pole} to be $M=(V,E,S_{1},\ldots,S_{n})$ where
$S_{1},\ldots,S_{n}$ is a partition of the set $S(M)$ into pairwise disjoint
sets $S_{i}$ such that $\left\vert S_{i}\right\vert =k_{i}$ for $i=1,\ldots,n$
and $k_{1}+\cdots+k_{n}=k.$ Each set $S_{i}$ is called a \emph{connector} of
$M,$ or more precisely $k_{i}$\emph{-connector}.

\begin{definition}
A \emph{supervertex} (resp. \emph{superedge}) is a multipole with three (resp.
two) connectors.
\end{definition}

We define two particular supervertices $A$ and $A^{\prime}$ as in Figure
\ref{Fig_superVert}. Let $G$ be a snark and $u_{1},u_{2}$ a pair of
non-adjacent vertices of $G.$ A superedge $G_{u_{1},u_{2}}$ is obtained from
$G$ by removing vertices $u_{1}$ and $u_{2},$ and then replacing all edges
incident to a vertex $u_{i}$ in $G$ by semiedges in $G_{u_{1},u_{2}}$ which
form a connector $S_{i}$, for $i=1,2.$

\begin{definition}
A \emph{proper} superedge is either an isolated edge or a multipole
$G_{u_{1},u_{2}}$ where $G$ is a snark.
\end{definition}

\noindent A proper superedge can be defined in a wider sense (see
\cite{KocholSuperposition}), but for the purposes of the present paper this
simple definition suffices. Now, let $G=(V,E)$ be a cubic graph. Replace each
edge $e\in E$ by a superedge $\mathcal{E}(e)$ and each vertex $v\in V$ by a
supervertex $\mathcal{V}(v)$ so that the following holds: if $v\in V$ is
incident to $e\in E$, then semiedges of one connector of $\mathcal{V}(v)$ are
identified with semiedges of one connector from $\mathcal{E}(e)$ and these two
connectors must have the same cardinality. Notice that $\mathcal{V}(v)$ (resp.
$\mathcal{E}(e)$) does not have to be the same for every vertex $v$ (resp.
edge $e),$ so $\mathcal{V}$ (resp. $\mathcal{E}$) is a mapping of the set of
vertices $V(G)$ to a set of supervertices (resp. superedges).

This procedure results in a cubic graph which is called a \emph{superposition}
of $G$ with $\mathcal{V}$ and $\mathcal{E}$ and denoted by $G(\mathcal{V}%
,\mathcal{E})$. Moreover, if $\mathcal{E}(e)$ is proper for every $e\in E$,
then $G(\mathcal{V},\mathcal{E})$ is called a \emph{proper} superposition of
$G$. Since a multipole consisting of a single vertex and three semiedges is a
supervertex, and an isolated edge is proper superedge, some of the vertices
and edges of $G$ may be superposed by themselves. The following theorem is
established in \cite{KocholSuperposition}. There, it is stated for snarks of
girth $\geq5$ which are cyclically $4$-edge-connected, but it also holds for
snarks of smaller girth and smaller edge-connectivity.

\begin{theorem}
\label{Tm_Kochol}If $G$ is a snark and $G(\mathcal{V},\mathcal{E})$ is a
proper superposition of $G$, then $G(\mathcal{V},\mathcal{E})$ is a snark.
\end{theorem}

In this paper, we will establish that Conjecture \ref{Con_normal} holds for
snarks $G(\mathcal{V},\mathcal{E})$ obtained for particular supervertices
$\mathcal{V}$ and superedges $\mathcal{E}$ superposed on vertices and edges of
a cycle in $G.$

\paragraph{Normal colorings of superpositioned snarks.}

Normal colorings of a family of superpositioned snarks are considered in
\cite{SedSkrePaper1}. There, for any snark $G$ and any even cycle $C$ in $G,$
a superposition $G_{C}(\mathcal{A},\mathcal{\mathcal{B)}}$ is considered. For
every vertex $u_{i}$ of $C,$ a supervertex $\mathcal{A}(u_{i})$ is either $A$
or $A^{\prime}$ from Figure \ref{Fig_superVert}. For every edge $e_{i}$ of
$C$, a superedge $\mathcal{B}(e_{i})$ is the supervertex $(P_{10})_{u,v}$. All
other vertices and edges of $G$ are superpositioned by themselves. For many
such superpositions the existence of a normal $5$-edge-coloring is
established, and the approach naturally extends to superposition of any even
subgraph of $G$ consisting of even length cycles. The approach of
\cite{SedSkrePaper1} could not be extended to odd cycles, since $P_{10}$ is a
small snark. Thus, in this paper we consider a similar superposition with
Flower snarks taken as superedges instead of the Petersen graph. This class of
superpositions by Flower snarks has already been considered in \cite{Kochol2,
KocholFlower1, kineziEuropean, MacajovaSkovieraHypohamiltonian}, only in a
different context.

Flower snarks are first introduced by Isaacs in \cite{Isaacs1975}, and in
\cite{Hagglund2014} it is established that Flower snarks have a normal
$5$-edge-coloring. Flower snarks are used as superedges in \cite{Kochol2,
KocholFlower1} in order to obtain cyclically $6$-edge-connected snarks. There,
a $6$-cycle of the Petersen graph is superpositioned by supervertices $A$ or
$A^{\prime}$ and by superedges obtained from Flower snarks. A generalized
construction, where any snark $G$ is taken instead of the Petersen graph and
any cycle $C$ instead of a $6$-cycle, is used in
\cite{MacajovaSkovieraHypohamiltonian} to obtain hypohamiltonian snarks for
almost all orders. In \cite{kineziEuropean} it is established that for such
superposition the Berge-Fulkerson Conjecture holds.

In this paper, we first consider a superposition obtained from any snark $G$
by superpositioning vertices of a cycle $C$ in $G$ by supervertices $A$ and
$A^{\prime}$, edges of $C$ by superedges $H_{x,y}$, where $H$ is any snark and
not necessarily the same for every edge of $C,$ and $x,y$ is any pair of
nonadjacent vertices of $H.$ For such superpositions, two sufficient
conditions for the existence of a normal $5$-edge-coloring are given. The
first condition is applied to snarks superpositioned by hypohamiltonian
snarks, which yields a normal $5$-edge-coloring for many of them. To be more
precise, the condition yields a normal $5$-edge-coloring for all
hypohamiltonian snarks used as superedges, but only for some of the possible
ways of connecting them. In particular, since Flower snarks are
hypohamiltonian, this consequently yields a normal $5$-edge-coloring for many
(but not all) snarks superpositioned by Flower snarks. The second sufficient
condition is more demanding, but its application yields a normal
$5$-edge-colorings for all superpositions by Flower snarks.

Finally, since the existence of a Petersen coloring implies the existence of
the six perfect matchings required by the Berge-Fulkerson Conjecture, the
results of this paper imply the results from \cite{kineziEuropean}.

\section{Preliminaries}

Let $G$ be a snark, and $C=u_{0}u_{1}\cdots u_{g-1}u_{0}$ a cycle of length
$g$ in $G.$ Denote the edges of a cycle $C$ by $e_{i}=u_{i}u_{i+1}$ for
$i=0,\ldots,g-1,$ where indices are taken modulo $g$. Also, let $v_{i}$ denote
the neighbor of $u_{i}$ distinct from $u_{i-1}$ and $u_{i+1},$ and let
$f_{i}=u_{i}v_{i}.$

\begin{figure}[h]
\begin{center}
\includegraphics[scale=0.6]{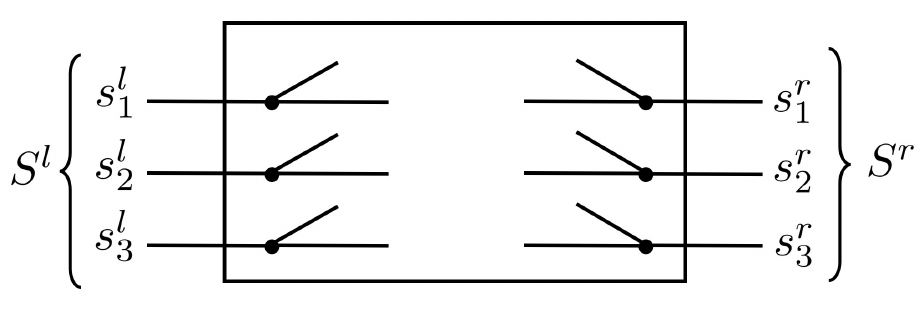}
\end{center}
\caption{A schematic depiction of a superedge $\mathcal{B}_{i}$ of a
superposition $G_{C}(\mathcal{A},\mathcal{B}),$ with its left and right
connector and their semiedges. Throughout the paper, we will assume such
depiction.}%
\label{Fig_scehamticConnector}%
\end{figure}

Recall that we defined two particular supervertices $A$ and $A^{\prime}$ as in
Figure \ref{Fig_superVert}. For superedges we take $H_{x,y},$ where $H$ is a
snark and $x,y$ a pair of nonadjacent vertices of $H.$ Notice that $H_{x,y}$
has a pair of $3$-connectors, the connectors $S_{x}$ and $S_{y}$ containing
the three semiedges which are halves of the three edges of $H$ incident to $x$
and $y$, respectively. Let $G_{C}(\mathcal{A},\mathcal{B})$ be a superposition
of $G$ such that $\mathcal{A}(u_{i})\in\{A,A^{\prime}\}$, $\mathcal{B}%
(e_{i})\in\{H_{x,y}:H$ is a snark and $x,y$ a pair of nonadjacent vertices of
$H\}$ for every $i=0,\ldots,g-1,$ and all other vertices and edges of $G$ are
superpositioned by themselves. The family of such superpositions is denoted by
$\mathcal{G}_{C}(\mathcal{A},\mathcal{B}).$ We will write shortly
$\mathcal{A}_{i}$ for $\mathcal{A}(u_{i})$ and $\mathcal{B}_{i}$ for
$\mathcal{B}(e_{i}).$ Notice that in general we do not obtain a same graph
$G_{C}(\mathcal{A},\mathcal{B}),$ up to isomorphism, if the connector $S_{x}$
of $\mathcal{B}_{i}$ is identified with $\mathcal{A}_{i}$ and $S_{y}$ with
$\mathcal{A}_{i+1},$ as when $S_{x}$ is connected to $\mathcal{A}_{i+1}$ and
$S_{y}$ to $\mathcal{A}_{i}.$ In other words, we may consider a superedge
$\mathcal{B}_{i}$ as being directed from $S_{x}$ to $S_{y}.$ To be more
precise, a connector of $\mathcal{B}_{i}$ which is identified with a connector
of $\mathcal{A}_{i}$ (resp. $\mathcal{A}_{i+1}$) will be called the
\emph{left} (resp. \emph{right}) connector and denoted by $S^{l}$ (resp.
$S^{r}$), its three semiedges will be called the \emph{left} (resp.
\emph{right}) semiedges and denoted by $s_{1}^{l},$ $s_{2}^{l},$ $s_{3}^{l}$
(resp. $s_{1}^{r},$ $s_{2}^{r},$ $s_{3}^{r}$). A scheme of superedge which
will be used throughout the paper is shown in Figure
\ref{Fig_scehamticConnector}. For a superedge $\mathcal{B}_{i}$ obtained from
a snark $H$ by removing nonadjacent vertices $x$ and $y,$ we write
$\mathcal{B}_{i}=H_{x,y}$ (resp. $\mathcal{B}_{i}=H_{y,x}$) in case when
$S_{x}$ (resp. $S_{y}$) is the left connector.

\begin{figure}[h]
\begin{center}
\includegraphics[scale=0.6]{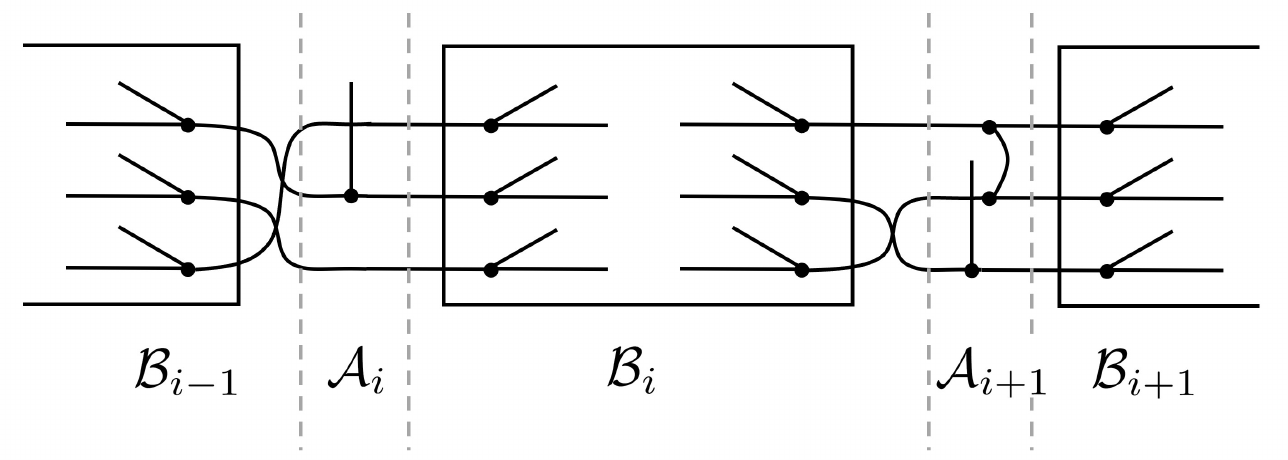}
\end{center}
\caption{The figure shows how a superedge $\mathcal{B}_{i}$ is connected to
$\mathcal{B}_{i-1}$ and $\mathcal{B}_{i+1}$ given a permutation $p_{i}$ and a
dock $d_{i}$ associated with $\mathcal{B}_{i}.$ In this example,
$p_{i}=(1,3,2)$ and $d_{i}=2.$ From the figure it is also visible that
$p_{i-1}=(2,3,1)$ and $d_{i+1}=3.$}%
\label{Fig_permutationConnection}%
\end{figure}

It is convenient to assume that the connection of a supervertex $\mathcal{A}%
_{i}$ with superedges $\mathcal{B}_{i-1}$ and $\mathcal{B}_{i}$ arises so that
first the right semiedges of $\mathcal{B}_{i-1}$ are identified with the left
semiedges of $\mathcal{B}_{i}$ and then one of the arising edges is subdivided
to obtain the vertex $u_{i}$ of $\mathcal{A}_{i}.$ Since the semiedge
identification of $\mathcal{B}_{i-1}$ and $\mathcal{B}_{i}$ can be done in
various ways, as shown in Figure \ref{Fig_permutationConnection}, we will
associate a permutation $p_{i-1}$ of the set $\{1,2,3\}$ with a superedge
$\mathcal{B}_{i-1}$, so that the right superedges $s_{1}^{r},$ $s_{2}^{r},$
$s_{3}^{r}$ of $\mathcal{B}_{i-1}$ are first permuted according to $p_{i-1}$
and then identified with left semiedges $s_{1}^{l},$ $s_{2}^{l},$ $s_{3}^{l}$
of $\mathcal{B}_{i},$ i.e. a right semiedge $s_{p_{i-1}^{-1}(j)}^{r}$ of
$\mathcal{B}_{i-1}$ is identified with a left $s_{j}^{l}$ of $\mathcal{B}%
_{i}.$ A permutation $p_{i-1}$ will be called a \emph{semiedge permutation}.
For example, in Figure \ref{Fig_permutationConnection} it holds that
$p_{i-1}=(2,3,1)$ and $p_{i}=(1,3,2)$. When $p_{i-1}$ is clear from the
context, we will for brevity sake denote $p_{i-1}^{-1}(j)$ by $j^{-}.$

Assuming that $s_{j^{-}}^{r}s_{j}^{l}$ denotes an edge arising by semiedge
identification, one $j\in\{1,2,3\}$ has to be chosen as the index of the edge
which is to be subdivided by $u_{i}.$ Since $j$ is the index of left semiedges
in $\mathcal{B}_{i},$ we will denote that choice by $j=d_{i}$ and associate it
with $\mathcal{B}_{i}.$ The choice $d_{i}$ of a left semiedge of
$\mathcal{B}_{i}$ to be connected to the vertex $u_{i}$ of $\mathcal{A}_{i}$
is called a \emph{dock} index, and the left semiedge $s_{d_{i}}^{l}$ is called
a \emph{dock} semiedge. For example, in Figure \ref{Fig_permutationConnection}
it holds that $d_{i}=2$ and $d_{i+1}=3.$

To conclude, each superedge $\mathcal{B}_{i}$ is associated with a permutation
$p_{i}$ which determines how the right semiedges of $\mathcal{B}_{i}$ connect
to left semiedges of $\mathcal{B}_{i+1},$ and a dock index $d_{i}$ which
determines which of the left semiedges of $\mathcal{B}_{i}$ is connected to a
vertex $u_{i}$ of $\mathcal{A}_{i}.$

\paragraph{Submultipoles and their compatible colorings.}

We will next consider normal $5$-edge-colorings of a superposition
$G_{C}(\mathcal{A},\mathcal{B}),$ for which purpose we first need to define
the notions of submultipole and the restriction of a coloring to a
submultipole. Let $M=(V,E,S)$ be a multipole. A multipole $M^{\prime
}=(V^{\prime},E^{\prime},S^{\prime})$ is a \emph{submultipole} of $M$ if
$V^{\prime}\subseteq V,$ $E^{\prime}\subseteq E$ and $S^{\prime}\subseteq
S\cup E_{S}$ where $E_{S}$ is the set consisting of halves of edges from $E.$
For a subset $V^{\prime}\subseteq V,$ we say that a multipole $M^{\prime
}=M[V^{\prime}]$ is an \emph{induced} submultipole of $M$ if $M^{\prime}$ has
a set of vertices $V^{\prime},$ set of edges $E^{\prime}$ consisting of all
edges $e$ of $M$ with both end-vertices in $V^{\prime},$ and set of semiedges
$S^{\prime}$ which consists of all semiedges of $M$ whose only end-vertex
belongs to $V^{\prime}$ and halves of all edges of $E$ which have precisely
one end-vertex in $V^{\prime}.$

Further, let $\sigma$ be a normal $5$-edge-coloring of a cubic multipole $M.$
Then the \emph{restriction} $\sigma^{\prime}=\left.  \sigma\right\vert
_{M^{\prime}}$ of $\sigma$ to $M^{\prime}$ is defined by $\sigma^{\prime
}(e)=\sigma(e)$ for $e\in E^{\prime},$ $\sigma^{\prime}(s)=\sigma(s)$ for
$s\in S^{\prime}\cap S$ and $\sigma^{\prime}(s)=\sigma(e_{s})$ for $s\in
S^{\prime}\backslash S$ where $e_{s}$ denotes the edge of $M$ such that $s$ is
a semiedge of $e_{s}.$ Let $M_{1},\ldots,M_{k}$ be cubic submultipoles of a
cubic multipole $M,$ let $\sigma_{i}$ be a normal $5$-edge-coloring of a
submultipole $M_{i}$ for $i=1,\ldots,k,$ and let $M^{\prime}$ be a
submultipole of $M$ induced by $\cup_{i=1}^{k}V(M_{i}).$ Colorings $\sigma
_{i}$ are said to be \emph{compatible}, if there exists a normal
$5$-edge-coloring $\sigma^{\prime}$ of $M^{\prime}$ such that $\left.
\sigma^{\prime}\right\vert _{M_{i}}=\sigma_{i}$ for every $i=1,\ldots,k.$

\begin{figure}[h]
\begin{center}%
\begin{tabular}
[t]{llll}%
a) & \raisebox{-0.9\height}{\includegraphics[scale=0.6]{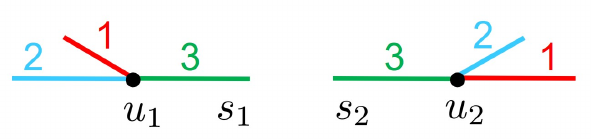}} & b) &
\raisebox{-0.9\height}{\includegraphics[scale=0.6]{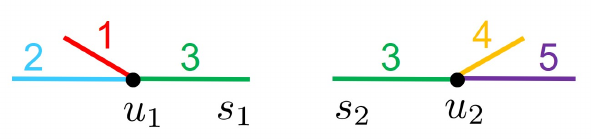}}\\
c) & \raisebox{-0.9\height}{\includegraphics[scale=0.6]{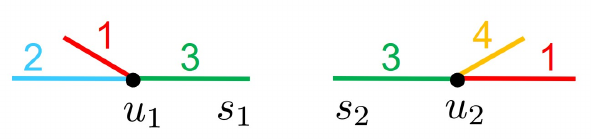}} & d) &
\raisebox{-0.9\height}{\includegraphics[scale=0.6]{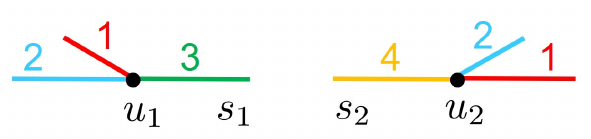}}
\end{tabular}
\end{center}
\caption{All four figures show the same color scheme $\sigma_{1}[s_{1}]$ of a
semiedge $s_{1}$. As for the color scheme $\sigma_{2}[s_{2}],$ it holds that:
a) $\sigma_{2}[s_{2}]$ is consistent with $\sigma_{1}[s_{1}]$ since they are
the same, b) $\sigma_{2}[s_{2}]$ is consistent with $\sigma_{1}[s_{1}]$ since
they coincide on $s_{1}$ and $s_{2}$ and have complementary pair of incident
colors, c) $\sigma_{2}[s_{2}]$ is not consistent with $\sigma_{1}[s_{1}]$
since they have neither the same nor complementary incident colors, d)
$\sigma_{2}[s_{2}]$ is not consistent with $\sigma_{1}[s_{1}]$ since they do
not coincide on $s_{1}$ and $s_{2}.$}%
\label{Fig_compatible_semiedges}%
\end{figure}

We introduce the notion of color schemes of semiedges, illustrated by Figure
\ref{Fig_compatible_semiedges}. Let $M$ be a cubic multipole, and $\sigma$ a
normal $5$-edge-coloring of $M.$ For a semiedge $s\in S(M)$ incident to a
vertex $v_{s}\in V(M),$ a\emph{ color scheme} $\sigma\lbrack s]$ is defined by
$\sigma\lbrack s]=(i,\{j,k\}),$ where $i=\sigma(s)$ and $j,k$ are the two
remaining colors incident to the vertex $v_{s}$ in $\sigma.$

\begin{definition}
Let $s_{1}$ (resp. $s_{2}$) be a semiedge of a cubic multipole $M_{1}$ (resp.
$M_{2}$) which has normal $5$-edge-coloring $\sigma_{1}$ (resp. $\sigma_{2}$).
Color schemes
\[
\sigma_{1}[s_{1}]=(i_{1},\{j_{1},k_{1}\})\text{\quad and\quad}\sigma_{2}%
[s_{2}]=(i_{2},\{j_{2},k_{2}\})
\]
are \emph{consistent} if $i_{1}=i_{2}$ and $\{j_{1},k_{1}\}$ is equal either
to $\{j_{2},k_{2}\}\ $or\ to $\{1,2,3,4,5\}\backslash\{i_{2},j_{2},k_{2}\}.$
If $\sigma_{1}[s_{1}]$ and $\sigma_{2}[s_{2}]$ are consistent, we write
$\sigma_{1}[s_{1}]\approx\sigma_{2}[s_{2}].$
\end{definition}

Now, let $\tilde{G}\in\mathcal{G}_{C}(\mathcal{A},\mathcal{B})$ be a
superposition of $G$, $\mathcal{B}_{i}$ a superedge of $\tilde{G},$ and
$\tilde{\sigma}_{i}$ a $5$-edge-coloring of $\mathcal{B}_{i}.$ A\emph{ color
scheme} of the left and right connector of $\mathcal{B}_{i}$ is defined by
\[
\tilde{\sigma}_{i}[S^{l}]=(\tilde{\sigma}_{i}[s_{1}^{l}],\tilde{\sigma}%
_{i}[s_{2}^{l}],\tilde{\sigma}_{i}[s_{3}^{l}])\text{\quad and\quad}%
\tilde{\sigma}_{i}[S^{r}]=(\tilde{\sigma}_{i}[s_{1}^{r}],\tilde{\sigma}%
_{i}[s_{2}^{r}],\tilde{\sigma}_{i}[s_{3}^{r}]),
\]
respectively. A \emph{color scheme} of a superedge $\mathcal{B}_{i}$ is
defined by
\[
\tilde{\sigma}_{i}[\mathcal{B}_{i}]=(\tilde{\sigma}_{i}[S^{l}],\tilde{\sigma
}_{i}[S^{r}]).
\]
Let $\mathcal{B}_{i}$ be a superedge of $G_{C}(\mathcal{A},\mathcal{B}),$ and
$\tilde{\sigma}_{i}$ and $\tilde{\sigma}_{i}^{\prime}$ a pair of normal
$5$-edge-colorings of $\mathcal{B}_{i}.$ Colorings $\tilde{\sigma}_{i}$ and
$\tilde{\sigma}_{i}^{\prime}$ are \emph{consistent} on the left connector
$S^{l}$ of $\mathcal{B}_{i}$, which is denoted by $\tilde{\sigma}_{i}%
[S^{l}]\approx\tilde{\sigma}_{i}^{\prime}[S^{l}]$, if $\tilde{\sigma}%
_{i}[s_{j}^{l}]\approx\tilde{\sigma}_{i}^{\prime}[s_{j}^{l}]$ for every
$j=1,2,3.$ The consistency of $\tilde{\sigma}_{i}$ and $\tilde{\sigma}%
_{i}^{\prime}$ on the right connector $S^{r}$ of $\mathcal{B}_{i}$ is defined
and denoted analogously. Finally, colorings $\tilde{\sigma}_{i}$ and
$\tilde{\sigma}_{i}^{\prime}$ are \emph{consistent} on $\mathcal{B}_{i}$,
which is denoted by $\tilde{\sigma}_{i}[\mathcal{B}_{i}]\approx\tilde{\sigma
}_{i}^{\prime}[\mathcal{B}_{i}],$ if $\tilde{\sigma}_{i}$ and $\tilde{\sigma
}_{i}^{\prime}$ are consistent on both right and left connector of
$\mathcal{B}_{i}.$ If $\tilde{\sigma}_{i}$ and $\tilde{\sigma}_{i}^{\prime}$
are consistent on $\mathcal{B}_{i},$ we also say that color schemes
$\tilde{\sigma}_{i}[\mathcal{B}_{i}]$ and $\tilde{\sigma}_{i}^{\prime
}[\mathcal{B}_{i}]$ are consistent.

\section{Sufficient conditions for the existence of a normal $5$-edge-coloring
of a superposition}

We define a submultipole $M_{\mathrm{int}}$ of $\tilde{G}\in\mathcal{G}%
_{C}(\mathcal{A},\mathcal{B})$ as a submultipole induced by the set of
vertices $V(G)\backslash V(C)\subseteq V(\tilde{G}).$ Also, every
$\mathcal{B}_{i}$ is obviously a submultipole of $\tilde{G}.$ For a given
normal $5$-edge-coloring $\tilde{\sigma}$ of $\tilde{G},$ the corresponding
colorings $\tilde{\sigma}_{\mathrm{int}}$ of $M_{\mathrm{int}}$ and
$\tilde{\sigma}_{i}$ of $\mathcal{B}_{i}$ for all $i$ are defined as the
restrictions of $\tilde{\sigma}$ to a corresponding multipole, i.e.
$\tilde{\sigma}_{\mathrm{int}}=\left.  \tilde{\sigma}\right\vert
_{M_{\mathrm{int}}}$ and $\tilde{\sigma}_{i}=\left.  \tilde{\sigma}\right\vert
_{\mathcal{B}_{i}}.$ Notice that $\tilde{\sigma}_{\mathrm{int}}$ and
$\tilde{\sigma}_{i}$ are normal $5$-edge-colorings, and they are compatible in
the sense that there exists a normal $5$-edge-coloring of $\tilde{G}$ (that
being $\tilde{\sigma},$ of course) such that its restriction to the
corresponding multipoles is equal to $\tilde{\sigma}_{\mathrm{int}}$ and
$\tilde{\sigma}_{i}.$ Moreover, $\tilde{\sigma}_{\mathrm{int}}$ and all
$\tilde{\sigma}_{0},\ldots,\tilde{\sigma}_{g-1}$ uniquely identify
$\tilde{\sigma}.$ In our considerations, we will go in the reverse direction,
we will construct colorings $\tilde{\sigma}_{\mathrm{int}}$ and $\tilde
{\sigma}_{i}$ for all $i,$ and then establish that they are compatible and
thus we will obtain a desired normal $5$-edge-coloring of $\tilde{G}.$

We first introduce the so called right coloring, in which all three right
semiedges have a same color scheme, i.e. they are mutually consistent. Also,
when a normal $5$-edge-coloring of a superposition cannot be constructed using
only right colorings, a pair of semiedges $\mathcal{B}_{i-1}$ and
$\mathcal{B}_{i}$ will be considered together, where $\mathcal{B}_{i}$ will be
colored by a so called left coloring, so the names "left" and "right" for
colorings correspond to the position of $\mathcal{B}_{i-1}$ and $\mathcal{B}%
_{i}$ in illustrations.

\begin{figure}[h]
\begin{center}%
\begin{tabular}
[c]{lll}%
a) & \hspace*{-0.3cm}b) & \hspace*{-0.3cm}c)\\
\includegraphics[scale=0.55]{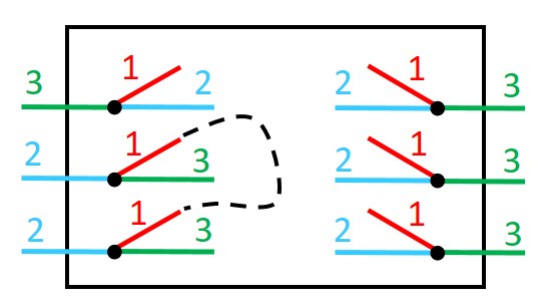} & \hspace*{-0.3cm}%
\includegraphics[scale=0.55]{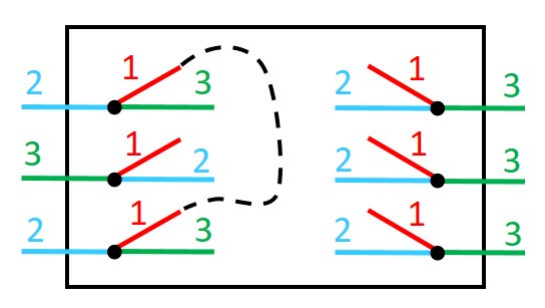} & \hspace*{-0.3cm}%
\includegraphics[scale=0.55]{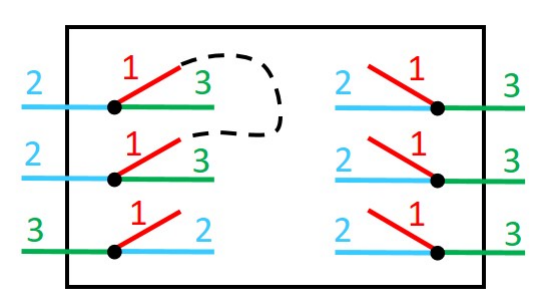}
\end{tabular}
\end{center}
\caption{The figure shows a color scheme $\kappa_{j}$ for: a) $j=1,$ b) $j=2$,
c) $j=3.$ Also, the dashed curve represents the Kempe chain $P^{l}.$\ }%
\label{Fig_schemesKappa}%
\end{figure}

\subsection{Right colorings and the first sufficient condition}

Let $\kappa_{j}$, for $j\in\{1,2,3\},$ denote the superedge color scheme from
Figure \ref{Fig_schemesKappa}. A $j$\emph{-right coloring} $R_{j}%
^{i}((2,1),3)$ of a superedge $\mathcal{B}_{i}$ is any normal $5$%
-edge-coloring of $\mathcal{B}_{i}$ consistent with a color scheme $\kappa
_{j}$ in which there exists a Kempe $(2,1)$-chain $P^{l}$ connecting the pair
of left semiedges distinct from $s_{j}^{l}.$ The Kempe chain $P^{l}$ is also
illustrated by Figure \ref{Fig_schemesKappa}. A superedge $\mathcal{B}_{i}$ is
said to be $j$\emph{-right} if it has a $j$-right coloring $R_{j}%
^{i}((2,1),3)$. In this notation $R_{j}^{i}((2,1),3),$ the number $3$ denotes
the color of all three right semiedges and the left semiedge $s_{j}^{l},$ and
$(2,1)$ denotes the colors of Kempe chain $P^{l}$ starting with color $2$
which is the color of the remaining two left semiedges distinct from
$s_{j}^{l}.$ Notice that in this kind of coloring all three right semiedges
are pairwise consistent. The word "right" is added the prefix "$j$-", since
precisely one of the three left semiedges, that being semiedge $s_{j}^{l},$ is
consistent with them.

A superedge $\mathcal{B}_{i}$ is called:\ 

\begin{itemize}
\item \emph{dock-right}, if it is a $j$-right for $j=d_{i},$

\item \emph{doubly-right,} if it is $j$-right for at least two distinct values
of $j,$

\item \emph{fully-right,} if it is $j$-right for every $j\in\{1,2,3\}.$
\end{itemize}

A coloring of $\mathcal{B}_{i}$ which can be obtained from $R_{j}%
^{i}((2,1),3)$ by replacing colors $(2,1)$ along the Kempe chain $P^{l}$ with
$(t_{1},t_{2})\in\{(1,2),(4,5),(5,4)\}$ will be also called a $j$\emph{-right
coloring}, and it will be denoted by $R_{j}^{i}((t_{1},t_{2}),3)$. Notice that
$R_{j}^{i}((t_{1},t_{2}),3)$ remains normal and it is consistent with
$R_{j}^{i}((2,1),3)$ on the right connector of the superedge, i.e. all three
right semiedges remain pairwise consistent in all $j$-right colorings of
$\mathcal{B}_{i}.$ Finally, notice that $j$-right coloring is a $5$%
-edge-coloring even though color schemes from Figure \ref{Fig_schemesKappa}
contain only three colors. So, any coloring of $\mathcal{B}_{i}$ which can be
obtained from a $j$-right coloring by a color permutation $c$, which is a
permutation of the set of $5$ colors, will also be called a $j$\emph{-right
}coloring. Such a coloring is denoted by $R_{j}^{i}((c(t_{1}),c(t_{2}%
)),c(3)).$

\begin{figure}[h]
\begin{center}
\includegraphics[scale=0.55]{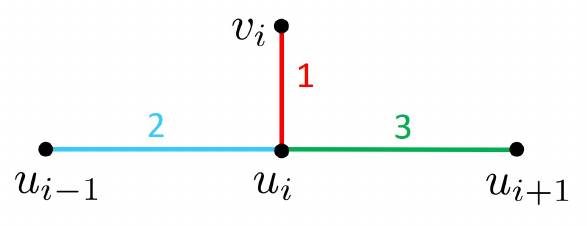}
\end{center}
\caption{A coloring $\sigma$ of the edges incident to a vertex $u_{i}$ of a
cycle $C$ in $G.$ For such a coloring $\sigma,$ all colorings of
$\mathcal{B}_{i}$ consistent with color schemes from Figure
\ref{Fig_schemesKappa} are $\sigma$-compatible.}%
\label{Fig_colSigma1}%
\end{figure}

Recall that edges of $G$ incident to vertices $u_{i}$ of the cycle $C$ in $G$
are denoted by $e_{i}=u_{i}u_{i+1}$ and $f_{i}=u_{i}v_{i}.$ Let $\sigma$ be a
normal $5$-edge-coloring of $G.$ A $j$-right coloring of $\mathcal{B}_{i}$ is
$\sigma$\emph{-compatible} if on the each of the three right semiedges as well
as on $s_{j}^{l}$ it is consistent with $(\sigma(e_{i}),\{\sigma(f_{i}%
),\sigma(e_{i-1})\})$, on the remaining two left semiedges it is consistent
with $(\sigma(e_{i-1}),\{\sigma(f_{i}),\sigma(e_{i})\})$ and there exists a
Kempe $(\sigma(e_{i-1}),\sigma(f_{i}))$-chain $P^{l}$ connecting these two
left semiedges, i.e. a coloring $R_{j}^{i}((\sigma(f_{i}),\sigma
(e_{i-1})),\sigma(e_{i})).$ A $j$-right coloring of $\mathcal{B}_{i}$ which is
$\sigma$-compatible will often be denoted only by $R_{j}^{i}.$ For example,
any $j$-right coloring consistent with color schemes shown in Figure
\ref{Fig_schemesKappa} is $\sigma$-compatible with $\sigma$ as in Figure
\ref{Fig_colSigma1}.

If, for a given $j$, a superedge $\mathcal{B}_{i}$ has a $j$-right coloring
which is $\sigma$-compatible for $\sigma$ as in Figure \ref{Fig_colSigma1},
then for any other $\sigma$ a corresponding $\sigma$-compatible $j$-right
coloring of $\mathcal{B}_{i}$ is obtained by color permutation. Consequently,
it is sufficient to check whether, for a given $j$, a superedge $\mathcal{B}%
_{i}$ has a $j$-right coloring in colors $1$, $2$,$\ 3$ as in Figure
\ref{Fig_schemesKappa}.

For a $\sigma$-compatible $j$-right coloring $R_{j}^{i}$ of $\mathcal{B}_{i},$
its \emph{complementary} $j$-right coloring $\bar{R}_{j}^{i}$ of
$\mathcal{B}_{i}$ is the $j$-right coloring which can be obtained from
$R_{j}^{i}$ by swapping colors $\sigma(e_{i-1})$ and $\sigma(f_{i})$ along the
Kempe chain $P^{l},$ i.e.
\[
\bar{R}_{j}^{i}=R_{j}^{i}((\sigma(f_{i}),\sigma(e_{i-1})),\sigma(e_{i})).
\]
We are now in a position to establish the first sufficient condition for the
existence of a normal $5$-edge-coloring of a superposition $\tilde{G}%
\in\mathcal{G}_{C}(\mathcal{A},\mathcal{B}).$

\begin{figure}[h]
\begin{center}
\includegraphics[scale=0.55]{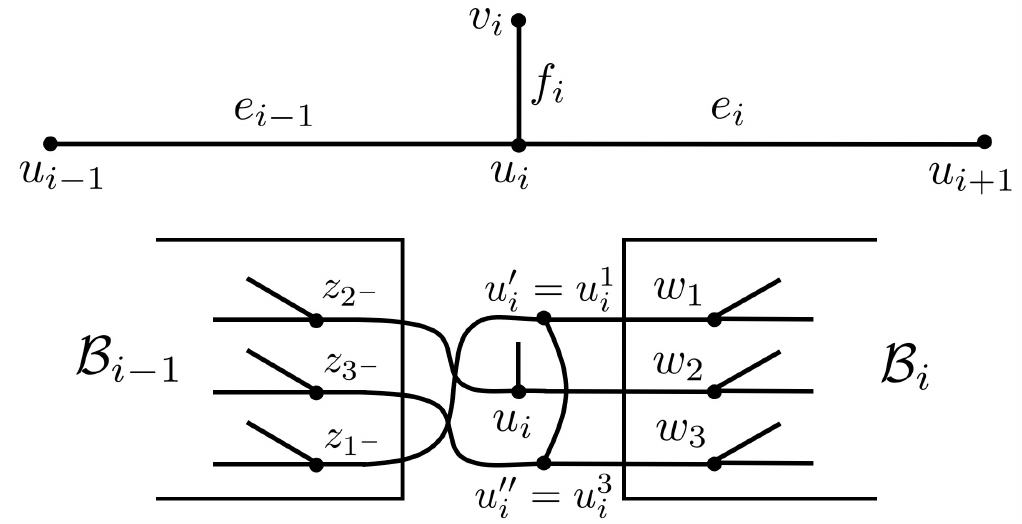}
\end{center}
\caption{The figure shows a supervertex $\mathcal{A}_{i},$ and superedges
$\mathcal{B}_{i-1}$ and $\mathcal{B}_{i}$ with $p_{i-1}=(2,3,1)$ and $d_{i}=2$
in a superposition $G_{C}(\mathcal{A},\mathcal{\mathcal{B}})$ and the
corresponding vertex notation.}%
\label{Fig_verticesWZ}%
\end{figure}

\begin{theorem}
\label{Tm_oneSE}Let $G$ be a snark with a normal $5$-edge-coloring $\sigma,$
$C$ a cycle of length $g$ in $G,$ and $\tilde{G}\in\mathcal{G}_{C}%
(\mathcal{A},\mathcal{B})$ a superposition of $G$. If the superedge
$\mathcal{B}_{i}$ is dock-right for every $i=0,\ldots,g-1$, then $\tilde{G}$
has a normal $5$-edge-coloring.
\end{theorem}

\begin{proof}
We will construct a normal $5$-edge-coloring $\tilde{\sigma}$ of $\tilde{G}$
by defining compatible colorings $\tilde{\sigma}_{\mathrm{int}}$ of
$M_{\mathrm{int}}$ and $\tilde{\sigma}_{i}$ of $\mathcal{B}_{i}$ for all $i.$
We first introduce some necessary notation. For a left semiedge $s_{j}^{l}$ of
$\mathcal{B}_{i},$ let $w_{j}$ denote the vertex of $\mathcal{B}_{i}$ incident
to it as illustrated in Figure \ref{Fig_verticesWZ}. Notice that $s_{j}^{l}$
of $\mathcal{B}_{i}$ is identified with a right semiedge  $s_{p_{i-1}^{-1}%
(j)}^{r}$ of $\mathcal{B}_{i-1}$ in the case of $\mathcal{A}_{i}=A$, or it is
a part of an incident edge in the case of $\mathcal{A}_{i}=A^{\prime}$. We
also introduce $j^{-}=p_{i-1}^{-1}(j)$. So, according to this simplified
notation, a right semiedge $s_{j^{-}}^{r}$ of $\mathcal{B}_{i-1}$ connects to
$s_{j}^{l}$ of $\mathcal{B}_{i}$. Further, let $z_{j^{-}}$ denote the vertex
of $\mathcal{B}_{i-1}$ to which $s_{j^{-}}^{r}$ is incident to, see Figure
\ref{Fig_verticesWZ}. Finally, recall that edges of $G$ incident to vertices
of $C$ are denoted by $e_{i}=u_{i}u_{i+1}$ and $f_{i}=u_{i}v_{i}.$

\bigskip\noindent\textbf{Claim 1.} \textit{If }$\mathcal{A}_{i}=A,$\textit{
then }$\tilde{\sigma}_{\mathrm{int}}=\left.  \sigma\right\vert
_{M_{\mathrm{int}}},$\textit{ }$\tilde{\sigma}_{i-1}=R_{j_{i-1}}^{i-1}%
$\textit{ and }$\tilde{\sigma}_{i}=R_{d_{i}}^{i}$\textit{ are compatible in
}$\tilde{G}$ \textit{for every }$j_{i-1}\in\{1,2,3\}.$

\medskip\noindent Let $M_{i}$ denote a submultipole of $\tilde{G}$ induced by
$V(M_{\mathrm{int}})\cup V(\mathcal{B}_{i-1})\cup V(\mathcal{A}_{i})\cup
V(\mathcal{B}_{i}).$ Assume that for an edge-coloring $\tilde{\sigma}$ of
$M_{i}$ it holds that $\left.  \tilde{\sigma}\right\vert _{\mathcal{B}_{i-1}%
}=\tilde{\sigma}_{i-1},$ $\left.  \tilde{\sigma}\right\vert _{\mathcal{B}_{i}%
}=\tilde{\sigma}_{i},$ and $\left.  \tilde{\sigma}\right\vert
_{M_{\mathrm{int}}}=\tilde{\sigma}_{\mathrm{int}}.$ We have to establish that
$\tilde{\sigma}$ is well defined and normal. Notice that $\tilde{\sigma
}(w_{d_{i}}u_{i})=\tilde{\sigma}_{i}(s_{d_{i}}^{l})=\sigma(e_{i})$,
$\tilde{\sigma}(z_{d_{i}^{-}}u_{i})=\tilde{\sigma}_{i-1}(s_{d_{i}^{-}}%
^{r})=\sigma(e_{i-1})$ and $\tilde{\sigma}(u_{i}v_{i})=\sigma(f_{i})$. We
obtain
\[
\tilde{\sigma}(u_{i})=\{\sigma(e_{i-1}),\sigma(e_{i}),\sigma(f_{i}%
)\}=\sigma(u_{i}),
\]
so $\tilde{\sigma}$ is proper at $u_{i}$ since $\sigma$ is proper at $u_{i}.$
Further, notice that for every $j\in\{1,2,3\}$ it holds that%
\[
\tilde{\sigma}_{i-1}(z_{j^{-}})=\{\sigma(e_{i-2}),\sigma(e_{i-1}%
),\sigma(f_{i-1})\}\text{\quad and\quad}\tilde{\sigma}_{i}(w_{j}%
)=\{\sigma(e_{i-1}),\sigma(e_{i}),\sigma(f_{i})\}.
\]
So, the edge $z_{d_{i}^{-}}u_{i}$ is normal by $\tilde{\sigma}$ since
$e_{i-1}$ is normal by $\sigma$. Also, the edge $u_{i}w_{d_{i}}$ is poor by
$\tilde{\sigma}$ since $\tilde{\sigma}(u_{i})=\tilde{\sigma}_{i}(w_{j}).$

For $j\not =d_{i}$ it holds that
\[
\tilde{\sigma}(z_{j^{-}}w_{j})=\tilde{\sigma}_{i-1}(s_{j^{-}}^{r}%
)=\tilde{\sigma}_{i}(s_{j}^{l})=\tilde{\sigma}(e_{i-1}).
\]
So, $\tilde{\sigma}$ is well defined. Given the values of $\tilde{\sigma
}_{i-1}(z_{j^{-}})$ and $\tilde{\sigma}_{i}(w_{j})$ noted above, it follows
that $z_{j^{-}}w_{j}$ is normal by $\tilde{\sigma}$ since $e_{i-1}$ is normal
by $\sigma.$ This establishes the claim.

\bigskip\noindent\textbf{Claim 2.} \textit{If }$\mathcal{A}_{i}=A^{\prime}%
,$\textit{ then }$\tilde{\sigma}_{\mathrm{int}}=\left.  \sigma\right\vert
_{M_{\mathrm{int}}},$\textit{ }$\tilde{\sigma}_{i-1}=R_{j_{i-1}}^{i-1}%
$\textit{ and }$\tilde{\sigma}_{i}=\bar{R}_{d_{i}}^{i}$\textit{ are compatible
in }$\tilde{G}$\textit{ for any }$j_{i-1}\in\{1,2,3\}.$

\medskip\noindent Let $M_{i}$ and $\tilde{\sigma}$ be defined as in Claim 1.
The properness of $\tilde{\sigma}$ at $u_{i}$ and the normality of the edges
of $\tilde{G}$ incident to $u_{i}$ is established in the same way as in Claim
1. We define $\tilde{\sigma}(u_{i}^{\prime}u_{i}^{\prime\prime})=\sigma
(e_{i}).$ For the sake of notation consistency, let us denote $u_{i}^{\prime}$
(resp. $u_{i}^{\prime\prime}$) also by $u_{i}^{j}$, where $j$ is the index of
the left semiedge $s_{j}^{l}$ of $\mathcal{B}_{i}$ which connects to
$u_{i}^{j},$ see Figure \ref{Fig_verticesWZ}. Thus, for $j\not =d_{i}$ we
have
\[
\tilde{\sigma}(u_{i}^{j}w_{j})=\tilde{\sigma}_{i}(s_{j}^{l})=\sigma
(f_{i})\text{\quad and\quad}\tilde{\sigma}(u_{i}^{j}z_{j^{-}})=\tilde{\sigma
}_{i}(s_{j^{-}}^{r})=\sigma(e_{i-1}).
\]
This implies $\tilde{\sigma}(u_{i}^{j})=\sigma(u_{i}),$ so $\tilde{\sigma}$ is
proper at $u_{i}^{j}$ since $\sigma$ is proper at $u_{i}.$ Moreover, from
$\tilde{\sigma}_{i}(w_{j})=\sigma(u_{i})=\tilde{\sigma}(u_{i}^{j})$ it follows
that the edge $w_{j}u_{i}^{j}$ is poor by $\tilde{\sigma}$. Given the values
of $\tilde{\sigma}_{i-1}(z_{j^{-}})$ and $\tilde{\sigma}(u_{i}^{j}),$ it also
follows that the edge $z_{j^{-}}u_{i}^{j}$ is normal by $\tilde{\sigma}$ since
$e_{i-1}$ is normal by $\sigma.$ This establishes Claim 2.

\bigskip Next, for every $i=0,\ldots,g-1$, we define
\[
\tilde{\sigma}_{i}=\left\{
\begin{array}
[c]{ll}%
R_{d_{i}}^{i} & \text{if }\mathcal{A}_{i}=A,\\
\bar{R}_{d_{i}}^{i} & \text{if }\mathcal{A}_{i}=A^{\prime}.
\end{array}
\right.
\]
Since $R_{d_{i}}^{i}$ and $\bar{R}_{d_{i}}^{i}$ are consistent on the right
connector $S^{r}$ of $\mathcal{B}_{i}$, Claims 1 and 2 imply that
$\tilde{\sigma}_{\mathrm{int}}$ and $\tilde{\sigma}_{i}$, for all $i$, are
compatible in $\tilde{G}$. This further implies that $\tilde{G}$ has a normal
$5$-edge-coloring, and we are finished.
\end{proof}

From the above theorem, we immediately have the next result as a corollary.

\begin{corollary}
\label{Cor_oneSE}Let $G$ be a snark with a normal $5$-edge-coloring $\sigma,$
$C$ a cycle of length $g$ in $G,$ and $\tilde{G}\in\mathcal{G}_{C}%
(\mathcal{A},\mathcal{B})$ a superposition of $G$. If the superedge
$\mathcal{B}_{i}$ is fully-right for every $i=0,\ldots,g-1$, then $\tilde{G}$
has a normal $5$-edge-coloring.
\end{corollary}

\begin{figure}[h]
\begin{center}%
\begin{tabular}
[c]{lll}%
\hspace*{-0.4cm}\includegraphics[scale=0.64]{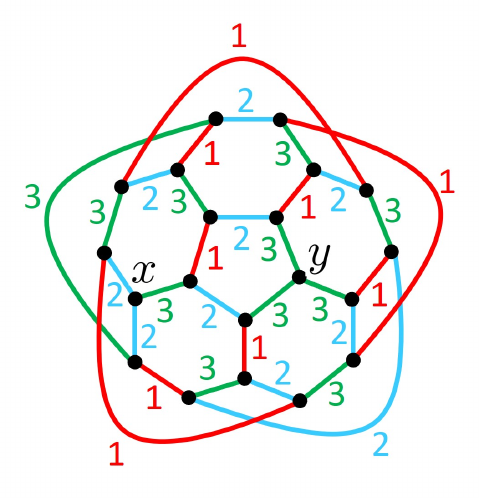} &
\includegraphics[scale=0.64]{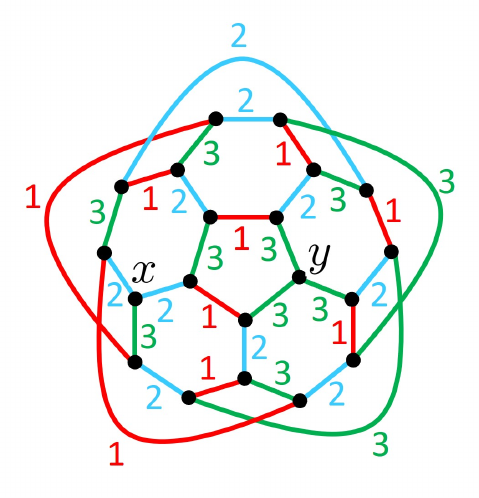} &
\includegraphics[scale=0.64]{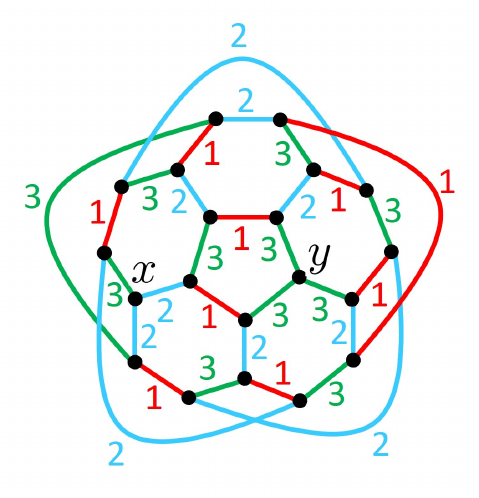}
\end{tabular}
\end{center}
\caption{The figure shows three distinct edge-colorings of the Flower snark
$H$ which are neither proper nor normal, but their restriction to the
superedge $H_{x,y}$ yields a right coloring where each of the left semiedges
is a dock in one of them. Thus, $H_{x,y}$ is fully-right.}%
\label{Fig_completelyRight}%
\end{figure}

\paragraph{Application to hypohamiltonian snarks.}

The sufficient condition of Theorem \ref{Tm_oneSE} applies to many snarks $H$
used as superedges. To see this, we first need to introduce the notion of
hypohamiltonian graphs. Namely, a graph $H$ is \emph{hypohamiltonian} if $H$
is not Hamiltonian and for every vertex $v\in V(H)$ it holds that $H-v$ is
Hamiltonian. In \cite{Fiorini1983}, it was established that Flower snarks are
hypohamiltonian, and in \cite{MacajovaSkovieraHypohamiltonian}%
\ hypohamiltonian snarks with cyclic connectivity 5 and 6 are constructed for
all but finitely many even orders. Thus, there are infinitely many
hypohamiltonian snarks, and for such snarks the following proposition holds.

\begin{proposition}
Let $H$ be a hypohamiltonian snark and $x,y$ a pair of nonadjacent vertices in
$H.$ Then $H_{x,y}$ is $j$-right for at least one $j\in\{1,2,3\}.$
\end{proposition}

\begin{proof}
Let $H^{\prime}=H-y.$ Since $H$ is hypohamiltonian, it follows that
$H^{\prime}$ is Hamiltonian. Let $C^{\prime}$ be a Hamiltonian cycle in
$H^{\prime}.$ Since $H$ is a cubic graph, it has an even number of vertices,
so $H^{\prime}$ has an odd number of vertices. Consequently, $C^{\prime}$ is
an odd length cycle in $H^{\prime}.$ Let $x_{1},$ $x_{2}$, and $x_{3}$ be the
three neighbors of $x$ in $H.$ Thus, $C^{\prime}$ contains precisely two out
of the three edges $x_{j}x,$ say edges $x_{1}x$ and $x_{2}x.$ Let
$H^{\prime\prime}$ (resp. $C^{\prime\prime}$) be obtained from $H^{\prime
}-x_{3}x$ (resp. $C^{\prime}$) by suppressing the vertex $x.$ Denote by
$e_{x}$ the edge of $H^{\prime\prime}$ obtained by suppressing $x.$

Now, notice that $H^{\prime\prime}$ is a graph of an even order and
$C^{\prime\prime}$ is a Hamiltonian cycle in it. Color the edges of
$C^{\prime\prime}$ in $H^{\prime\prime}$ alternately by $1$ and $2$ so that
$e_{x}$ is colored by $2,$ and all edges in $E(H^{\prime\prime})\backslash
E(C^{\prime\prime})$ are colored by $3$. This yields a proper $3$%
-edge-coloring $\sigma^{\prime\prime}$ of $H^{\prime\prime}.$ We construct the
corresponding coloring $\sigma$ of $H$ in the following way. All edges $e\in
E(H)\cap E(H^{\prime\prime})$ are colored in $H$ as in $H^{\prime\prime},$
edges $x_{1}x$ and $x_{2}x$ are colored by $2,$ and the three edges of $H$
incident to $y$ and $x_{3}x$ are colored by $3.$ Notice that $\sigma$ is
neither a proper nor a normal edge-coloring of $H$, but the restriction of
$\sigma$ to $H_{x,y}$ yields a $j$-right coloring of $H_{x,y}$ for $j=3$.
\end{proof}

In the above proposition, it is established that every hypohamiltonian snark
is $j$-right for at least one $j\in\{1,2,3\}$. Let $\tilde{G}$ be a
superposition of a snark $G$ such that every edge of a cycle $C$ is
superpositioned by a hypohamiltonian snark (which does not have to be the
same  one for all edges of $C$) in a way that the superedge obtained from the
hypohamiltonian snark is $j$-right for $j=d_{i}.$ For such a superposition,
Theorem \ref{Tm_oneSE} implies the existence of a normal $5$-edge-coloring of
the superposition $\tilde{G}$.

There also exist snarks $H$ and pairs of vertices $x,y$ in them, such that
$H_{x,y}$ is fully-right. For example, if one considers for $H$ the so-called
Flower snark $J_{5}$ shown in Figure \ref{Fig_completelyRight}, and for $x,y$
the pair of vertices designated in the figure, then it is evident from the
figure that such $H_{x,y}$ is fully-right. So, if $\tilde{G}=G_{C}%
(\mathcal{A},\mathcal{\mathcal{B}})$ is a superposition such that
$\mathcal{B}_{i}=H_{x,y}$ for every $i,$ then Corollary \ref{Cor_oneSE}
implies that $\tilde{G}$ has a normal $5$-edge-coloring.

\begin{figure}[h]
\begin{center}%
\begin{tabular}
[c]{lll}%
a) & \hspace*{-0.3cm}b) & \hspace*{-0.3cm}c)\\
\includegraphics[scale=0.55]{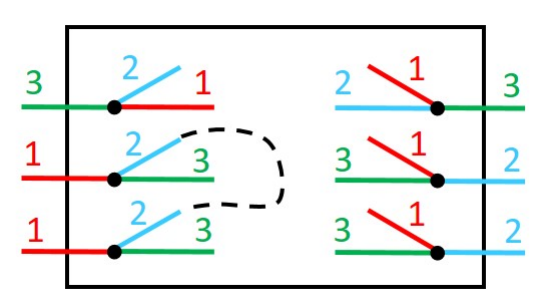} & \hspace*{-0.3cm}%
\includegraphics[scale=0.55]{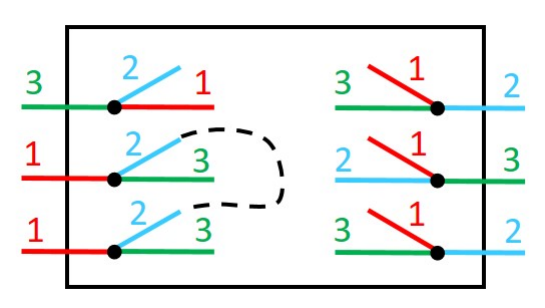} & \hspace*{-0.3cm}%
\includegraphics[scale=0.55]{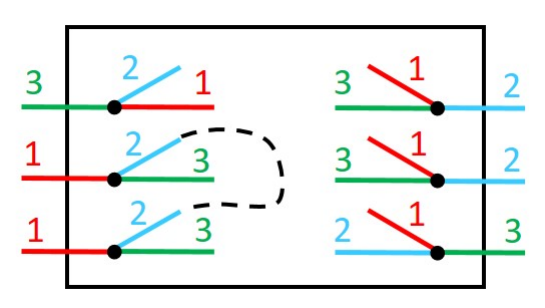}\\
d) & \hspace*{-0.3cm}e) & \hspace*{-0.3cm}f)\\
\includegraphics[scale=0.55]{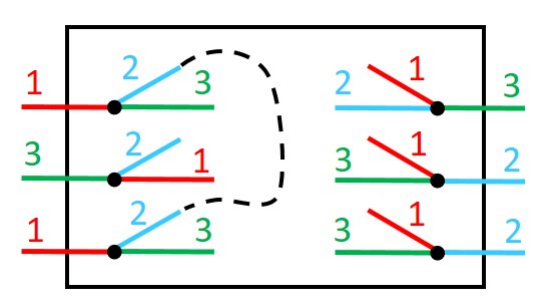} & \hspace*{-0.3cm}%
\includegraphics[scale=0.55]{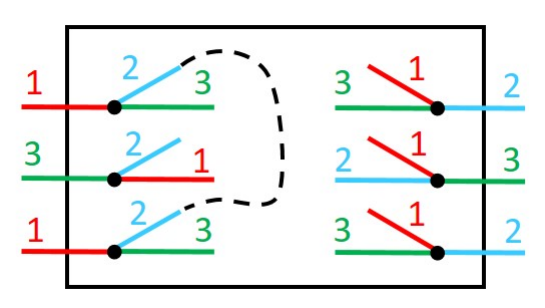} & \hspace*{-0.3cm}%
\includegraphics[scale=0.55]{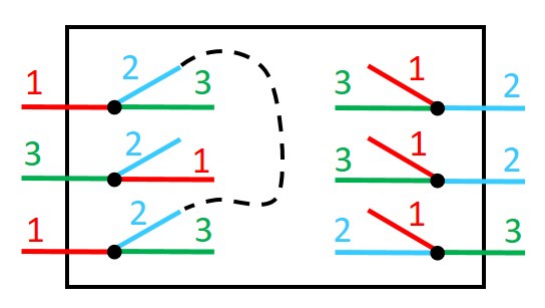}\\
g) & \hspace*{-0.3cm}h) & \hspace*{-0.3cm}i)\\
\includegraphics[scale=0.55]{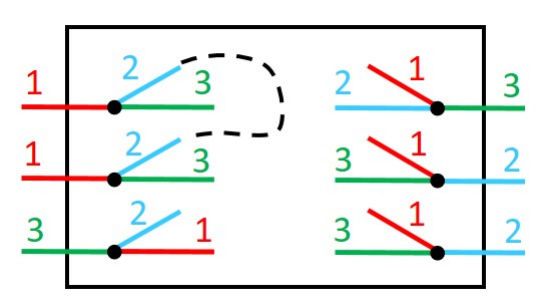} & \hspace*{-0.3cm}%
\includegraphics[scale=0.55]{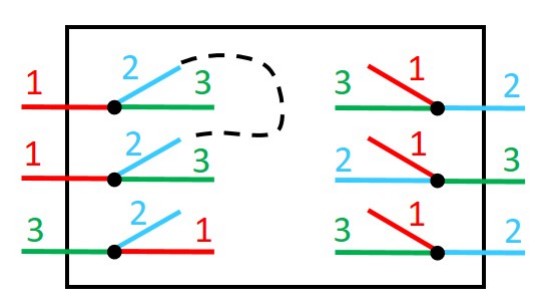} & \hspace*{-0.3cm}%
\includegraphics[scale=0.55]{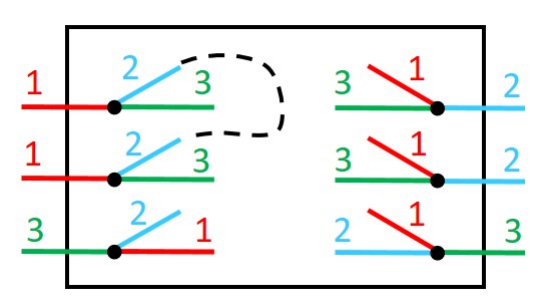}
\end{tabular}
\end{center}
\caption{The figure shows a color scheme $\tau_{j,k}$ for: a) $(j,k)=(1,1),$
b) $(j,k)=(1,2),$ c) $(j,k)=(1,3),$ d) $(j,k)=(2,1),$ e) $(j,k)=(2,2),$ f)
$(j,k)=(2,3),$ g) $(j,k)=(3,1),$ h) $(j,k)=(3,2),$ j) $(j,k)=(3,3).$}%
\label{Fig_schemeTau}%
\end{figure}

We will deal with Flower snarks in greater details in the next section. For
now, let us note that hypohamiltonian snarks are always dock-right  for at
least one dock out of three. Yet, they are often not fully-right.
Consequently, we need a refinement of the sufficient condition from Theorem
\ref{Tm_oneSE}. For that purpose, we introduce several more color schemes of a
superedge in the next section.

\subsection{Left colorings and the second sufficient condition}

Let $\tau_{j,k}$ denote a color scheme from Figure \ref{Fig_schemeTau}, for
$j,k\in\{1,2,3\}.$ Notice that in both the left and the right connectors of
$\tau_{j,k}$, two semiedges have the same color and the third semiedge is in a
distinct color. So, $j$ (resp. $k$) is the index of the left (resp. right)
semiedge in a color distinct from the color of the other two left (resp.
right) semiedges.

A $j,k$\emph{-left coloring} $L_{j,k}^{i}((1,2),3,(2,1))$ of $\mathcal{B}_{i}$
is any normal $5$-edge-coloring of $\mathcal{B}_{i}$ consistent with a color
scheme $\tau_{j,k}$ in which there exists a Kempe $(1,2)$-chain $P^{l}$
connecting the pair of left semiedges distinct from $s_{j}^{l}$. The chain
$P^{l}$ is also shown in Figure \ref{Fig_schemeTau}. Notice that the first
ordered pair $(1,2)$ of the notation $L_{j,k}^{i}((1,2),3,(2,1))$ indicates
that the two left semiedges distinct from $s_{j}^{l}$ are colored by $1$ and
that $P^{l}$ is a Kempe $(1,2)$-chain. The middle $3$ denotes the color of
semiedges $s_{j}^{l}$ and $s_{k}^{r}$.

The last ordered pair $(2,1)$ denotes that the remaining two right semiedges
distinct from $s_{k}^{r}$ are colored by $2$ and are incident to color $1$
beside $3$. In general, a left coloring is a $5$-edge-coloring. Yet, if
$L_{j,k}^{i}((1,2),3,(2,1))$ uses only $3$ colors, then the second ordered
pair being $(2,1)$ would imply the existence of a Kempe $(2,1)$-chain $P^{r}$
connecting a pair of right semiedges distinct from $s_{k}^{r}$, such that
$P^{l}$ and $P^{r}$ are vertex disjoint.

A superedge $\mathcal{B}_{i}$ is called \emph{doubly-left} if for every
$j\in\{1,2,3\}$ it holds that $\mathcal{B}_{i}$ has a left coloring
$L_{j,k}^{i}((1,2),3,(2,1))$ for at least two distinct values of $k$.

Again, we may assume that $L_{j,k}^{i}((1,2),3,(2,1))=L_{j,k}^{i}(T,3,(2,1))$
for $T=(1,2)$. Colors $T=(1,2)$ can be replaced along $P^{l}$ by any
$T^{\prime}$ from $\{(2,1),(4,5),(5,4)\}$. The coloring of $\mathcal{B}_{i}$
obtained in this way will also be called a $j,k$\emph{-left coloring} and
denoted by $L_{j,k}^{i}(T^{\prime},3,(2,1))$. Notice that $L_{j,k}%
^{i}(T^{\prime},3,(2,1))$ remains consistent with $L_{j,k}^{i}((1,2),3,(2,1))$
on the right connector for any $T^{\prime}.$ Finally, any coloring of
$\mathcal{B}_{i}$ which can be obtained from a $j,k$-left coloring by a color
permutation $c,$ which is a permutation of the set of five colors, will be
also called a $j,k$\emph{-left coloring}. Similarly as with right colorings,
if for a given $j$ and $k$ the superedge $\mathcal{B}_{i}$ has one $j,k$-left
coloring, then it has all $j,k$-left colorings.

\begin{figure}[h]
\begin{center}
\includegraphics[scale=0.55]{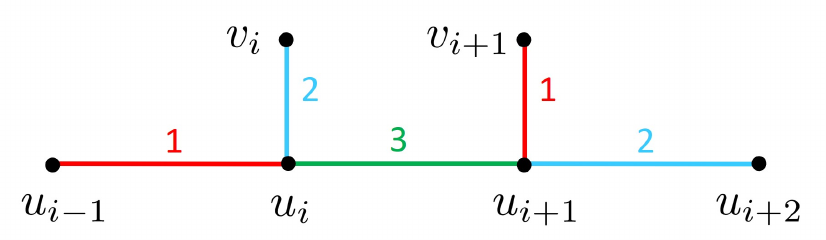}
\end{center}
\caption{The figure shows a coloring $\sigma$ of the edges incident to
vertices $u_{i}$ and $u_{i+1}$ of a cycle $C$ in $G$. All $j,k$-left colorings
of $\mathcal{B}_{i}$ consistent with color schemes from Figure
\ref{Fig_schemeTau} are $\sigma$-compatible with such a coloring $\sigma$ of
$G$.}%
\label{Fig_TauSigma}%
\end{figure}

Let $\sigma$ denote a normal $5$-edge-coloring of a snark $G$. A $j,k$-left
coloring of $\mathcal{B}_{i}$ is $\sigma$\emph{-compatible}\ if on $s_{j}^{l}$
and $s_{k}^{r}$ it is consistent with $(\sigma(e_{i}),\{\sigma(f_{i+1}%
),\sigma(e_{i+1})\}),$ on the two right semiedges distinct from $s_{k}^{r}$ it
is consistent with $(\sigma(e_{i+1}),\{\sigma(f_{i+1}),\sigma(e_{i})\}),$ and
on the two left semiedges distinct from $s_{j}^{l}$ it is consistent with
$(\sigma(e_{i-1}),\{\sigma(e_{i}),\sigma(f_{i})\})$ with a Kempe
$(\sigma(e_{i-1}),\sigma(f_{i}))$-chain $P^{l}$ connecting this pair of left
semiedges. In other words, it is $\sigma$-compatible if it is a coloring
$L_{j,k}^{i}(T_{1},\sigma(e_{i}),T_{2})$, where $T_{1}=(\sigma(e_{i-1}%
),\sigma(f_{i}))$ and $T_{2}=(\sigma(e_{i+1}),\sigma(f_{i+1})).$ Again, for a
$\sigma$-compatible $j,k$-left coloring we will often write just $L_{j,k}^{i}$
without denoting the colors. For example, any coloring consistent with color
schemes from Figure \ref{Fig_schemeTau} is $\sigma$-compatible for any
$\sigma$ as in Figure \ref{Fig_TauSigma}.

For a $\sigma$-compatible $j,k$-left coloring $L_{j,k}^{i}$ of $\mathcal{B}%
_{i}$, a \emph{complementary} $j,k$-left coloring $\bar{L}_{j,k}^{i}$ of
$\mathcal{B}_{i}$ is a $j,k$-left coloring which can be obtained from
$L_{j,k}^{i}$ by swapping colors $\sigma(e_{i-1})$ and $\sigma(f_{i-1})$ along
$P^{l}$, i.e.
\[
\bar{L}_{j,k}^{i}=L_{j,k}^{i}(T_{1}^{\prime},\sigma(e_{i}),T_{2})
\]
for $T_{1}^{\prime}=(\sigma(f_{i}),\sigma(e_{i-1}))$ and $T_{2}=(\sigma
(e_{i+1}),\sigma(f_{i+1})).$ Let us first make the following observation.

\begin{remark}
\label{Obs_leftRightConsistent}Colorings $L_{j,k}^{i}$ and $R_{j}^{i}$ are
consistent on the left connector of $\mathcal{B}_{i}$ for any $k\in\{1,2,3\}$.
The same holds for complementary colorings $\bar{L}_{j,k}^{i}$ and $\bar
{R}_{j}^{i}$.
\end{remark}

\begin{figure}[h]
\begin{center}
\includegraphics[scale=0.5]{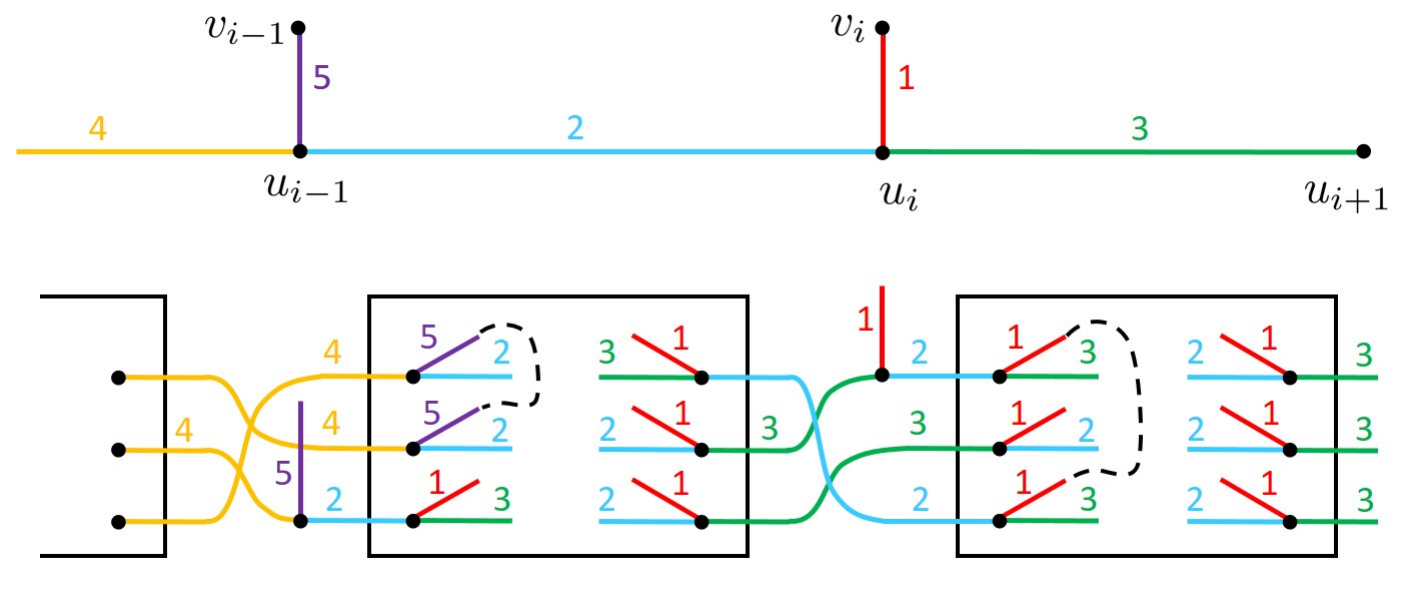}
\end{center}
\caption{Considering together $\mathcal{B}_{i-1}$ and $\mathcal{B}_{i}$ in
case when $\mathcal{B}_{i}$ is not dock-right. Superedge $\mathcal{B}_{i-1}$
is colored by a left coloring and $\mathcal{B}_{i}$ by a right coloring.}%
\label{Fig_doubleSE}%
\end{figure}

We will also need the following lemma.

\begin{lemma}
\label{Lema_leftRightCompatible}Let $G$ be a snark with a normal
$5$-edge-coloring $\sigma,$ $C$ a cycle of length $g$ in $G,$ and $\tilde
{G}\in\mathcal{G}_{C}(\mathcal{A},\mathcal{B})$ a superposition of $G$. Let
$\tilde{\sigma}_{i}$ be a coloring of $\mathcal{B}_{i}.$ If $\mathcal{A}%
_{i}=A$, then for $\{d_{i},j_{i},p_{i-1}(k_{i-1})\}=\{1,2,3\}$ and any
$j_{i-1}\in\{1,2,3\}$ it holds that $\tilde{\sigma}_{\mathrm{int}}=\left.
\sigma\right\vert _{M_{\mathrm{int}}},$ $\tilde{\sigma}_{i-1}=L_{j_{i-1}%
,k_{i-1}}^{i-1}$ and $\tilde{\sigma}_{i}=R_{j_{i}}^{i}$ are compatible in
$\tilde{G}.$
\end{lemma}

\begin{proof}
A pair of colorings $\tilde{\sigma}_{i-1}$ and $\tilde{\sigma}_{i}$ which
satisfy the conditions of the lemma is illustrated by Figure
\ref{Fig_doubleSE}. Let us assume that the vertices incident to the edges
containing the semiedges of $\mathcal{A}_{i}$ are denoted as in Figure
\ref{Fig_verticesWZ}. Namely, we denote indices of right semiedges in
$\mathcal{B}_{i-1}$ and their incident vertices by $j^{-}=p_{i-1}^{-1}(j).$
Assume that there exists a coloring $\tilde{\sigma}$ such that $\left.
\tilde{\sigma}\right\vert _{M_{\mathrm{int}}}=\tilde{\sigma}_{\mathrm{int}}$,
$\left.  \tilde{\sigma}\right\vert _{\mathcal{B}_{i-1}}=\tilde{\sigma}_{i-1}$
and $\left.  \tilde{\sigma}\right\vert _{\mathcal{B}_{i}}=\tilde{\sigma}_{i}$.
We have to establish that $\tilde{\sigma}$ is well defined and normal. Notice
that
\begin{equation}
\tilde{\sigma}_{i}(w_{j})=\tilde{\sigma}_{i-1}(z_{j^{-}})=\{\sigma
(e_{i-1}),\sigma(e_{i}),\sigma(f_{i})\}=\sigma(u_{i}) \label{For_wjzj}%
\end{equation}
for every $j\in\{1,2,3\}$. Also, notice that the index $j$, which is used to
denote the right semiedges of $\mathcal{B}_{i}$, takes its values from
$\{d_{i},j_{i},p_{i-1}(k_{i-1})\}$.

Let us first consider $j=d_{i}$, which implies $j\not =$ $p_{i-1}(k_{i-1})$,
i.e., $k_{i-1}\not =j^{-}$. In the left coloring $\tilde{\sigma}%
_{i-1}=L_{j_{i-1},k_{i-1}}^{i-1}$, it holds that $k_{i-1}\not =j^{-}$. By the
definition of the $\sigma$-compatible left coloring, we have $\tilde{\sigma
}_{i-1}(s_{j^{-}}^{r})=\sigma(e_{i})$. In the right coloring $\tilde{\sigma
}_{i}=R_{j_{i}}^{i}$, it holds that $j_{i}\not \in \{p_{i-1}(k_{i-1}),d_{i}%
\}$, i.e., $j_{i}\not =j$. By the definition of $\sigma$-compatible right
coloring, we have $\tilde{\sigma}_{i}(s_{j}^{l})=\sigma(e_{i-1})$. We conclude
that
\begin{align*}
\tilde{\sigma}(u_{i}w_{j})  &  =\tilde{\sigma}_{i}(s_{j}^{l})=\sigma
(e_{i-1}),\\
\tilde{\sigma}(u_{i}z_{j^{-}})  &  =\tilde{\sigma}_{i-1}(s_{j^{-}}^{r}
)=\sigma(e_{i}),\\
\tilde{\sigma}(u_{i}v_{i})  &  =\sigma(u_{i}v_{i})=\sigma(f_{i}),
\end{align*}
so we obtain $\tilde{\sigma}(u_{i})=\sigma(u_{i})$. Therefore, $\tilde{\sigma
}$ is proper at $u_{i}$ since $\sigma$ is proper at $u_{i}$. Combined with
(\ref{For_wjzj}), this also yields that edges $u_{i}w_{j}$, $u_{i}z_{j^{-}}$
are poor by $\tilde{\sigma}$, and $u_{i}v_{i}$ is normal by $\tilde{\sigma}$
since it is normal by $\sigma$.

Let us next consider $j=j_{i}\not \in \{p_{i-1}(k_{i-1}),d_{i}\}$. Since
$j\not =p_{i-1}(k_{i-1})$, it holds that $j^{-}\not =k_{i-1}$. So, by the
definition of $\sigma$-compatible left coloring, we have $\tilde{\sigma}%
_{i-1}(s_{j^{-}}^{r})=\sigma(e_{i})$. Also, since $j=j_{i}$, we have
$\tilde{\sigma}_{i}(s_{j}^{l})=\sigma(e_{i})$. We conclude
\[
\tilde{\sigma}(w_{j}z_{j^{-}})=\tilde{\sigma}_{i}(s_{j}^{l})=\tilde{\sigma
}_{i-1}(s_{j^{-}}^{r})=\sigma(e_{i}),
\]
which implies $\tilde{\sigma}$ is well defined on $w_{j}z_{j^{-}}$. Combined
with (\ref{For_wjzj}), it follows that $w_{j}z_{j^{-}}$ is poor by
$\tilde{\sigma}$.

We finally consider $j=p_{i-1}(k_{i-1})\not =j_{i}$. In this case,
$j=p_{i-1}(k_{i-1})$, which means $j^{-}=k_{i-1}$. By the definition of
$\sigma$-compatible left coloring, we have $\tilde{\sigma}_{i-1}(s_{j^{-}}%
^{r})=\sigma(e_{i-1})$. On the other hand, $j\not =j_{i}$. Therefore, by the
definition of $\sigma$-compatible right coloring, we have $\tilde{\sigma}%
_{i}(s_{j}^{l})=\sigma(e_{i-1})$. Thus, we have
\[
\tilde{\sigma}(w_{j}z_{j^{-}})=\tilde{\sigma}_{i}(s_{j}^{l})=\tilde{\sigma
}_{i-1}(s_{j^{-}}^{r})=\sigma(e_{i-1}),
\]
so $\tilde{\sigma}$ is well defined on $w_{j}z_{j^{-}}$. Combined with
(\ref{For_wjzj}), it follows that $w_{j}z_{j^{-}}$ is poor by $\tilde{\sigma}$.
\end{proof}

\begin{figure}[h]
\begin{center}
\includegraphics[scale=0.5]{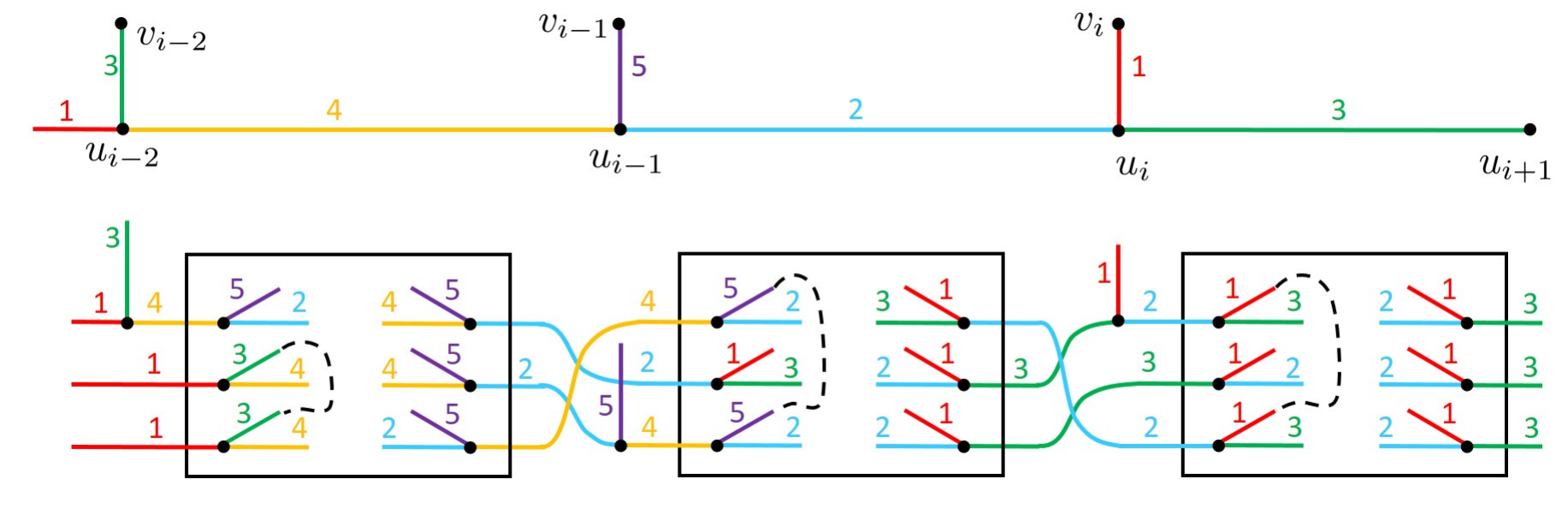}
\end{center}
\caption{Considering together $\mathcal{B}_{i-2}$, $\mathcal{B}_{i-1}$, and
$\mathcal{B}_{i}$, in the case when $\mathcal{B}_{i}$ is not dock-right.
Superedges $\mathcal{B}_{i-2}$ and $\mathcal{B}_{i-1}$ are colored by a left
coloring, and $\mathcal{B}_{i}$ by a right coloring.}%
\label{Fig_tripleSE}%
\end{figure}

We can now establish the following sufficient condition for the existence of a
normal $5$-edge-coloring of superpositions.

\begin{theorem}
\label{Tm_twoSE}Let $G$ be a snark with a normal $5$-edge-coloring $\sigma$,
$C$ a cycle of length $g$ in $G$, and $\tilde{G}\in\mathcal{G}_{C}
(\mathcal{A},\mathcal{B})$ a superposition of $G$. If the superedge
$\mathcal{B}_{i}$ is both doubly-right and doubly-left for every
$i=0,\ldots,g-1$, then $\tilde{G}$ has a normal $5$-edge-coloring.
\end{theorem}

\begin{proof}
Let us first explain the idea behind this proof, so that it is later easier to
follow. We will construct a normal $5$-edge-coloring $\tilde{\sigma}$ of
$\tilde{G}$ by defining a normal $5$-edge-coloring of larger "chunks" which
consist of two consecutive superedges. Namely, in the proof of Theorem
\ref{Tm_oneSE}, we considered single superedges and colored them by a right
coloring. We proved further that the left side of such a colored superedge is
compatible with the coloring of the previous superedge, and the right side is
compatible with the coloring of the subsequent superedge. Here, some
superedges might not have a desired right coloring. Each such superedge will
be considered together with the neighboring superedge, i.e. we will consider
two consecutive superedges as one larger "chunk". This "chunk" will be colored
so that the left side of the leftmost superedge in the "chunk" is colored the
same as the left side of a single superedge in Theorem \ref{Tm_oneSE}. The
right side of the rightmost superedge in a "chunk" will be colored the same as
the right side of a single superedge in Theorem \ref{Tm_oneSE}. Thus, these
"chunks" will be mutually compatible in the same way as single superedges are
in Theorem \ref{Tm_oneSE}.

Now we distinguish the following two cases with respect to whether $\tilde{G}$
has some superedge $\mathcal{A}_{i}$ that is $A^{\prime}$.

\bigskip\noindent\textbf{Case 1: }$\mathcal{A}_{i}=A$\emph{ for every
}$i=0,\ldots,g-1$. A desired coloring is obtained by the following algorithm.
We first explain the details of the three procedures used in the algorithm,
which are given input $i$, and afterwards, we describe the algorithm itself.
The three procedures are the following:

\bigskip\texttt{Procedure P1.} Consider the single dock-right superedge
$\mathcal{B}_{i}$, set $\tilde{\sigma}_{i}:=R_{d_{i}}^{i}$, and then set
$i:=i-1$.

\bigskip\texttt{Procedure P2.} Consider the pair of superedges $\mathcal{B}%
_{i-1}$ and $\mathcal{B}_{i}$ together, where $\mathcal{B}_{i}$ is not
dock-right, as illustrated by Figure \ref{Fig_doubleSE}. Since $\mathcal{B}%
_{i-1}$ is doubly-left, for $j_{i-1}=d_{i-1}$ there exists a $k_{i-1}%
\not =p_{i-1}^{-1}(d_{i})$ such that $\mathcal{B}_{i-1}$ has a left coloring
$L_{j_{i-1},k_{i-1}}^{i-1}$. Furthermore, $\mathcal{B}_{i}$ is doubly-right
but not dock-right. So, for $j_{i}\not \in \{p_{i-1}(k_{i-1}),d_{i}\}$, there
exists a right coloring $R_{j_{i}}^{i}$ of $\mathcal{B}_{i}$. Set%
\[
\tilde{\sigma}_{i-1}:=L_{j_{i-1},k_{i-1}}^{i-1}\text{\quad and\quad}%
\tilde{\sigma}_{i}:=R_{j_{i}}^{i}.
\]
Lemma \ref{Lema_leftRightCompatible} implies that $\tilde{\sigma
}_{\mathrm{int}}$, $\tilde{\sigma}_{i-1}$, and $\tilde{\sigma}_{i}$ are
compatible in $\tilde{G}$. After this is defined, set $i:=i-2$.

\bigskip\texttt{Procedure P3.} Consider the three superedges $\mathcal{B}%
_{i-2}$, $\mathcal{B}_{i-1}$, and $\mathcal{B}_{i}$ together, where
$\mathcal{B}_{i}$ is not dock-right, as illustrated by Figure
\ref{Fig_tripleSE}. Since $\mathcal{B}_{i-2}$ is doubly-left, for
$j_{i-2}=d_{i-2}$ there exists $k_{i-2} \not = p_{i-2}^{-1}(d_{i-1})$ such
that $\mathcal{B}_{i-2}$ has a left coloring $L_{j_{i-2},k_{i-2}}^{i-2}$.
Similarly, since $\mathcal{B}_{i-1}$ is doubly-left, for $j_{i-1}
\not \in \{p_{i-2}(k_{i-2}),d_{i-1}\}$, there exists $k_{i-1} \not =
p_{i-1}^{-1}(d_{i})$ such that $\mathcal{B}_{i-1}$ has a left coloring
$L_{j_{i-1},k_{i-1}}^{i-1}$. Finally, since $\mathcal{B}_{i}$ is doubly-right
but not dock-right, for the same $j_{i} \not \in \{p_{i-1}(d_{i-1}),d_{i}\}$,
there exists a right coloring $R_{j_{i}}^{i}$ of $\mathcal{B}_{i}$. Set
\[
\tilde{\sigma}_{i-2}:=L_{j_{i-2},k_{i-2}}^{i-2},\text{ }\tilde{\sigma}%
_{i-1}:=L_{j_{i-1},k_{i-1}}^{i-1}\text{ and }\tilde{\sigma}_{i}:=R_{j_{i}}%
^{i}.
\]
Lemma \ref{Lema_leftRightCompatible} implies that $\tilde{\sigma}_{i-1}$,
$\tilde{\sigma}_{i}$, and $\tilde{\sigma}_{\mathrm{int}}$ are compatible.
Together with Observation \ref{Obs_leftRightConsistent}, it also implies that
$\tilde{\sigma}_{i-2}$, $\tilde{\sigma}_{i-1}$, and $\tilde{\sigma
}_{\mathrm{int}}$ are compatible. After this is defined, set $i:=i-3$.

\bigskip

Let us now describe the flow of the algorithm. The algorithm starts with
$i:=g-1$. Suppose first that no superedge $\mathcal{B}_{i}$ of $\tilde{G}$ is
dock-right. If $g$ is even, apply \texttt{Procedure P2} repeatedly until
$i=-1$ is reached, and then the algorithm ends. If $g$ is odd, apply
\texttt{Procedure P2} until $i=2$. Then apply \texttt{Procedure P3} once for
the pair $\mathcal{B}_{1}$ and $\mathcal{B}_{0}$, which gives $i=-1$ at the end.

Suppose now that there exists a superedge $\mathcal{B}_{s}$ in $\tilde{G}$
which is dock-right. Without loss of generality, we may assume $s=0$. While $i
\geq3$, if $\mathcal{B}_{i}$ is dock-right, apply \texttt{Procedure P1}. If
$\mathcal{B}_{i}$ is not dock-right, apply \texttt{Procedure P2}. Since $i$ is
decreased by at most $2$ in \texttt{Procedures P1} and \texttt{P2}, the
execution will stop with the output $i \in\{1,2\}$. The algorithm concludes in
the following way:

\begin{itemize}
\item If $i=2$ and $\mathcal{B}_{2}$ is not dock-right, apply
\texttt{Procedure P3} once.

\item If $i=2$ and both $\mathcal{B}_{2}$ and $\mathcal{B}_{1}$ are
dock-right, apply \texttt{Procedure P1} three times.

\item If $i=2$ and $\mathcal{B}_{2}$ is dock-right and $\mathcal{B}_{1}$ is
not dock-right, apply first \texttt{Procedure P1} and then \texttt{Procedure
P2.}

\item If $i=1$ and $\mathcal{B}_{1}$ is dock-right, apply \texttt{Procedure
P1} twice.

\item If $i=1$ and $\mathcal{B}_{1}$ is not dock-right, apply
\texttt{Procedure P2} once.
\end{itemize}

\noindent In each of the above possibilities, the algorithm ends with $i=-1$.
This means that the algorithm has constructed $\tilde{\sigma}_{i}$ for every
$i=0,\ldots,g-1$. Let $Z$ denote the set of indices $i$ which were inputs for
\texttt{Procedures P1, P2} and \texttt{P3} in the execution of the algorithm.

It remains to establish that $\tilde{\sigma}_{i}$ and $\tilde{\sigma}_{i+1}$
are compatible for every $i\in Z$. Notice that $\tilde{\sigma}_{i}=R_{j}^{i}$
for some $j\in\{1,2,3\}$, and also $\tilde{\sigma}_{i+1}=R_{d_{i+1}}^{i+1}$ or
$\tilde{\sigma}_{i+1}=L_{d_{i+1},k}^{i+1}$ for some $k\in\{1,2,3\}$. In either
case, Observation \ref{Obs_leftRightConsistent} implies that $\tilde{\sigma
}_{i+1}$ is on the left side of $\mathcal{B}_{i+1}$, consistent with
$R_{d_{i+1}}^{i+1}$. Thus, $\tilde{\sigma}_{i}$ and $\tilde{\sigma}_{i+1}$ are
compatible by Claim 1 of Theorem \ref{Tm_oneSE}.

\bigskip\noindent\textbf{Case 2: }$\mathcal{A}_{i}=A^{\prime}$\emph{ for some
}$i\in\{0,\ldots,g-1\}$. Let $I$ denote the set of all indices $i$ such that
$\mathcal{A}_{i}=A^{\prime}$ in $\tilde{G}$, hence $I\not =\emptyset$. Denote
by $\tilde{G}^{A}$ a corresponding superposition obtained from $\tilde{G}$ by
setting $\mathcal{A}_{i}=A$ for every $i\in I$. Let $\tilde{\sigma}^{A}$ be a
normal $5$-edge-coloring of $\tilde{G}^{A}$ constructed as in Case 1. Denote
by $K$ the set of all indices $i$ such that $i-1$ is an output of Steps 1, 2,
or 3 in the execution of the algorithm of Case 1.

Notice that indices $i$ of $K$ represent the leftmost superedge $\mathcal{B}%
_{i}$ of each "chunk" of superedges considered together in the algorithm.
Thus, for every $i\in K$ from the procedures of the algorithm, it follows that
$\tilde{\sigma}_{i}^{A}=R_{d_{i}}^{i}$ or $\tilde{\sigma}_{i}^{A}%
=L_{d_{i},k_{i}}^{i}$. This means the index $j_{i}$ of the corresponding right
or left coloring of $\mathcal{B}_{i}$ is equal to the dock $d_{i}$. This also
implies that the two non-dock left semiedges of $\mathcal{B}_{i}$ are
connected by a Kempe chain $P^{l}$. Recall from Figure \ref{Fig_superVert}
that a supervertex $\mathcal{A}_{i}=A$ has two isolated edges, which connect
to the pair of non-dock left semiedges of $\mathcal{B}_{i}$. In the case of
$i\in K$, it also holds that these two isolated edges are colored by the same
color in $\tilde{\sigma}^{A}$, as illustrated by $\mathcal{A}_{i-1}$ in Figure
\ref{Fig_doubleSE} and $\mathcal{A}_{i-2}$ in Figure \ref{Fig_tripleSE}.

Let $T=\{0,\ldots,g-1\}\backslash K,$ i.e. $T$ denotes the set of all indices
not included in $K.$ Notice that for every $i\in T,$ the two isolated edges of
$\mathcal{A}_{i}=A$ are colored by different colors in $\tilde{\sigma}^{A},$
see $\mathcal{A}_{i}$ in Figure \ref{Fig_doubleSE} and $\mathcal{A}_{i-1},$
$\mathcal{A}_{i}$ in Figure \ref{Fig_tripleSE} for illustration. Also, it
holds that the Kempe chain $P^{l}$ of $\mathcal{B}_{i}$ does not connect the
two non-dock vertices. Instead, it connects a dock semiedge with one of the
remaining non-dock semiedge. This means that in the case of $i\in T,$ the
Kempe chain will not be of use to us, and instead we use the notion of a twist
defined in what follows.

Let $j_{1}$ and $j_{2}$ be the two distinct indices such that
$\{1,2,3\}\backslash\{d_{i}\}=\{j_{1},j_{2}\}.$ Notice that $\tilde{G}^{A}$
contains edges $w_{j_{1}}z_{j_{1}^{-}}$ and $w_{j_{2}}z_{j_{2}^{-}}.$ If, in
a  superposition $\tilde{G}^{A}$, the edges $w_{j_{1}}z_{j_{1}^{-}}$ and
$w_{j_{2}}z_{j_{2}^{-}}$ are replaced by $w_{j_{1}}z_{j_{2}^{-}}$ and
$w_{j_{2}}z_{j_{1}^{-}},$ this is called \emph{twisting} $\mathcal{A}_{i}$ in
$\tilde{G}^{A}.$ For a subset of indices $X \subseteq\{0,1,\ldots,g-1\}$, we
define the graph $\tilde{G}^{A,X}$ as the superposition obtained from
$\tilde{G}^{A}$ by twisting $\mathcal{A}_{i}$ for every $i \in X.$ Let
$\mathcal{G}^{A}$ denote the family of superpositions obtained from $\tilde
{G}^{A}$, which is defined by
\[
\mathcal{G}^{A}=\{\tilde{G}^{A,X}:X\subseteq\{0,1,\ldots,g-1\}\}.
\]
Since twisting does not influence the dock edge of any $\mathcal{B}_{i}$, the
algorithm of Case 1 will have the same set $K$ for every $\tilde{G}^{A,X}
\in\mathcal{G}^{A}$. Therefore, the "chunks" of superedges considered together
will  be the same for each $\tilde{G}^{A,X}$ from $\mathcal{G}^{A}$.

Let $\tilde{G}^{A,T}\in\mathcal{G}^{\prime}$ be the superposition obtained
from $\tilde{G}^{A}$ by twisting $\mathcal{A}_{i}$ for every $i\in T.$ Let
$\tilde{\sigma}^{A,T}$ be a normal $5$-edge-coloring of $\tilde{G}^{A,T}$
obtained by the algorithm from Case 1. Recall that $I$ denotes the set of all
indices $i$ such that $\mathcal{A}_{i}=A^{\prime}$ in $\tilde{G}.$ Let us
denote $I^{K}=I\cap K$ and $I^{T}=I\cap T.$ We now use $\tilde{\sigma}^{A,T}$
to construct a normal $5$-edge-coloring $\tilde{\sigma}$ of $\tilde{G}$ in the
following way. First, we define $\tilde{\sigma}_{i}$ of $\mathcal{B}_{i}$ by
slightly modifying $\tilde{\sigma}^{A,T}$, where it is important that in the
modification the color scheme of the right connector of $\mathcal{B}_{i}$
remains the same. Then we establish the compatibility of the consecutive
$\tilde{\sigma}_{i}$'s with $\tilde{\sigma}_{\mathrm{int}}.$ In order to do
this, we distinguish the following three cases.

\medskip\noindent\textbf{Subcase 2.a: }$i\in I^{K}.$ In this case,
$\tilde{\sigma}_{i}^{A,T}=R_{d_{i}}^{i}$ or $\tilde{\sigma}_{i}^{A,T}%
=L_{d_{i},k}^{i}$, i.e., $j_{i}=d_{i}$. The two isolated edges of
$\mathcal{A}_{i}=A$ within $\tilde{G}^{A,T}$ are colored with the same color
in $\tilde{\sigma}_{i}^{A,T}$. Consequently, the corresponding pair of left
semiedges of $\mathcal{B}_{i}$ are connected by a Kempe chain $P^{l}.$ This is
illustrated by the coloring of $\mathcal{A}_{i-1}$ from Figure
\ref{Fig_doubleSE}, since there it holds that $i-1\in K$. We define a coloring
$\tilde{\sigma}_{i}$ of $\mathcal{B}_{i}$ within $\tilde{G}$ as a coloring
obtained from $\tilde{\sigma}_{i}^{A,T}$ by swapping colors along $P^{l}$. The
edge $u_{i}^{\prime}u_{i}^{\prime\prime}$ within $\tilde{G}$ also has to be
colored. For that purpose, notice that the existence of Kempe chain $P^{l}$
connecting a pair of vertices $w_{j},$ for $j\not =d_{i},$ ensures that
$\tilde{\sigma}_{i}^{A,T}(w_{j})$ is the same for every $j\not =d_{i}$. We may
assume $\tilde{\sigma}_{i}^{A,T}(w_{j})=\{t_{1},t_{2},t_{3}\}$, where $t_{1}$
and $t_{2}$ are the colors along $P_{l}$. Thus, we define $\tilde{\sigma}%
_{i}(u_{i}^{\prime}u_{i}^{\prime\prime})=t_{3}$. Notice that  the color scheme
of the right connector of $\mathcal{B}_{i}$ remains the same  in
$\tilde{\sigma}_{i}$ as in $\tilde{\sigma}_{i}^{A,T}.$

Let us establish that $\tilde{\sigma}_{i}$, $\tilde{\sigma}_{i-1},$ and
$\tilde{\sigma}_{\mathrm{int}}=\left.  \sigma\right\vert _{M_{\mathrm{int}}}$
are compatible in $\tilde{G},$ assuming that $\tilde{\sigma}_{i-1}$ is
consistent with $\tilde{\sigma}_{i-1}^{A,T}$ on the right connector. For the
sake of notation consistency, $u_{i}^{\prime}$ and $u_{i}^{\prime\prime}$ are
denoted by $u_{i}^{j}$ if they are adjacent to $w_{j}$ of $\mathcal{B}_{i}.$
It is sufficient to establish properness at $u_{i}^{j}$ and the normality of
all edges incident to them. We may assume $\tilde{\sigma}_{i}(s_{j}^{l}
)=t_{1}$ for each $j\not =d_{i}$. Thus, $\tilde{\sigma}_{i}(w_{j}
)=\tilde{\sigma}_{i}(u_{i}^{j})=\{t_{1},t_{2},t_{3}\}$. Therefore,
$\tilde{\sigma}$ is proper at $u_{i}^{j}$ and edges $w_{j}u_{i}^{j}$ and
$u_{i}^{\prime} u_{i}^{\prime\prime}$ are poor in $\tilde{\sigma}$. This also
implies that an  edge $z_{j-}u_{i}^{j}$ has the same incident colors in
$\tilde{\sigma}$ as the  edge $z_{j^{-}}w_{j}$ in $\tilde{\sigma}^{A,T}$,
which implies normality.

\medskip\noindent\textbf{Subcase 2.b: }$i\in I^{T}.$ Notice that in this case
edges of $\mathcal{A}_{i}$ are twisted in $\tilde{G}^{A,T}$ compared to
$\tilde{G}^{A}$ (or $\tilde{G}$). Also, recall that in this case the two
isolated edges of $\mathcal{A}_{i}=A$ in $\tilde{G}^{A,T}$ are colored by two
different colors, which is illustrated by the coloring of $\mathcal{A}_{i}$
from Figure \ref{Fig_doubleSE} and also by the colorings of both
$\mathcal{A}_{i-1}$ and $\mathcal{A}_{i}$ from Figure \ref{Fig_tripleSE}. We
define the coloring $\tilde{\sigma}_{i}$ of $\mathcal{B}_{i}$ within
$\tilde{G}$ by $\tilde{\sigma}_{i}=\tilde{\sigma}_{i}^{A,T}.$ In order to
define the color of the edge $u_{i}^{\prime}u_{i}^{\prime\prime},$ let us
denote by $j_{1}$ and $j_{2}$ the two non-dock left indices of $\mathcal{B}%
_{i}$. Also, let us relabel $u_{i}^{\prime}$ and $u_{i}^{\prime\prime}$ by
$u_{i}^{j_{1}}$ and $u_{i}^{j_{2}}$ so that $u_{i}^{j_{l}}$ and $w_{j_{l}}$
are adjacent in $\tilde{G}$ (see Figure \ref{Fig_verticesWZ}). We may assume
that $\tilde{\sigma}_{i}(s_{d_{i}}^{l})=\tilde{\sigma}_{i}(s_{j_{1}}%
^{l})\not =\tilde{\sigma}_{i}(s_{j_{2}}^{l}).$ Since $\mathcal{B}_{i}$ is
colored by a $\sigma$-compatible coloring, it holds that $\tilde{\sigma}%
_{i}(s_{j_{2}}^{l})=\sigma(e_{i})$. This implies that $\tilde{\sigma}%
_{i}(w_{j_{2}})$ contains $\sigma(e_{i}).$ Let us denote $\tilde{\sigma}%
_{i}[s_{j_{2}}^{l}]\approx(\sigma(e_{i}),\{t_{1},t_{2}\}).$

\begin{figure}[h]
\begin{center}%
\begin{tabular}
[t]{ll}%
a) & \raisebox{-0.9\height}{\includegraphics[scale=0.55]{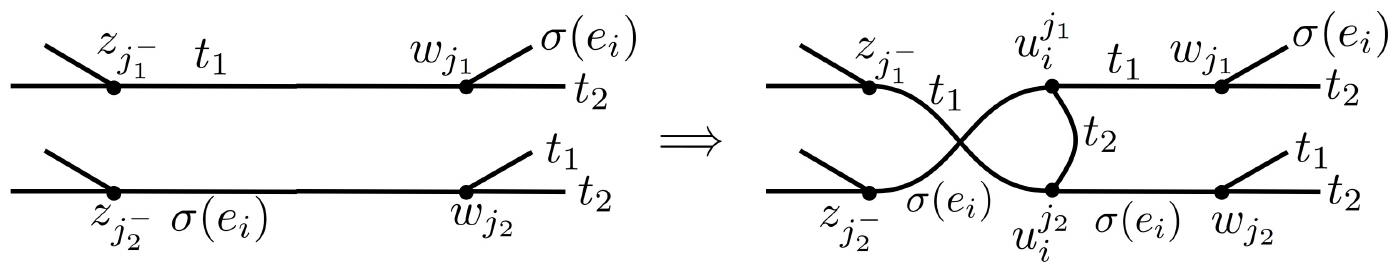}}\\
b) & \raisebox{-0.9\height}{\includegraphics[scale=0.55]{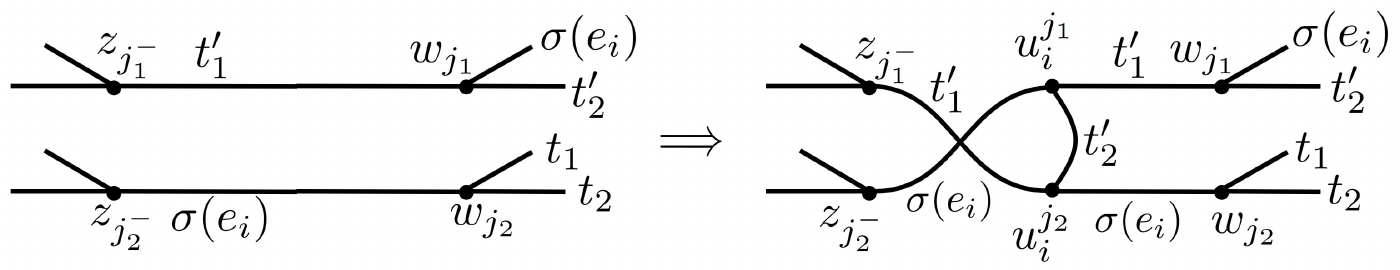}}
\end{tabular}
\end{center}
\caption{The figure shows the coloring $\tilde{\sigma}_{i}^{A,T}$ of the
superedge $\mathcal{A}_{i}=A$ in $\tilde{G}^{A,T}$ on the left, and the
coloring $\tilde{\sigma}_{i}$ of the twisted $\mathcal{A}_{i}=A^{\prime}$ in
$\tilde{G}$ on the right, where $\tilde{\sigma}_{i}[s_{j_{2}}^{l}%
]\approx(\sigma(e_{i}),\{t_{1},t_{2}\})$ and: a) $\tilde{\sigma}_{i}[s_{j_{1}%
}^{l}]\approx(t_{1},\{\sigma(e_{i}),t_{2}\})$, b) $\tilde{\sigma}_{i}%
[s_{j_{1}}^{l}]\approx(t_{1}^{\prime},\{\sigma(e_{i}),t_{2}^{\prime}\})$.}%
\label{Fig_twist}%
\end{figure}

Notice that by a right or a left color scheme from Figures
\ref{Fig_schemesKappa} and \ref{Fig_schemeTau}, respectively, it holds that
$\tilde{\sigma}_{i}[s_{j_{1}}^{l}]\approx(t_{1},\{\sigma(e_{i}),t_{2}\}).$
Moreover, in the case of $\tilde{\sigma}_{i}$ being a left coloring with
colors along $P^{l}$ replaced with the two colors not present in the scheme,
it holds that $\tilde{\sigma}_{i}[s_{j_{1}}^{l}]\approx(t_{1}^{\prime
},\{\sigma(e_{i}),t_{2}^{\prime}\})$, where $\{t_{1}^{\prime},t_{2}^{\prime
}\}=\{1,2,3,4,5\}\backslash\{\sigma(e_{i}),t_{1},t_{2}\}.$ Now, the value of
$\tilde{\sigma}_{i}$ on the edge $u_{i}^{j_{1}}u_{i}^{j_{2}}$ is defined
according to the following table:%
\[%
\begin{tabular}
[c]{|l|l||l|}\hline
$\tilde{\sigma}_{i}[s_{j_{2}}^{l}]$ & $\tilde{\sigma}_{i}[s_{j_{1}}^{l}]$ &
$\tilde{\sigma}_{i}(u_{i}^{j_{1}}u_{i}^{j_{2}})$\\\hline\hline
$(\sigma(e_{i}),\{t_{1},t_{2}\})$ & $(t_{1},\{\sigma(e_{i}),t_{2}\})$ &
$t_{2}$\\\cline{2-3}
& $(t_{1}^{\prime},\{\sigma(e_{i}),t_{2}^{\prime}\})$ & $t_{2}^{\prime}
$\\\hline
\end{tabular}
\ .
\]

\noindent These two cases are illustrated by Figure \ref{Fig_twist}. By
carefully checking all edges incident to $u_{i}^{j_{1}}$ and $u_{i}^{j_{2}}$
in each of the figures, one can verify that all of them are normal by
$\tilde{\sigma}_{i}$.

\medskip\noindent\textbf{Subcase 2.c: }$i\not \in I.$ Here, the coloring
$\tilde{\sigma}_{i}$ of $\mathcal{B}_{i}$ is defined simply as $\tilde{\sigma
}_{i}=\tilde{\sigma}_{i}^{A,T}$. In this subcase, $\tilde{\sigma}_{i}$,
$\tilde{\sigma}_{i-1}$, and $\tilde{\sigma}_{\mathrm{int}}$ are compatible in
$\tilde{G}$ since $\tilde{\sigma}_{i}^{A,T}$, $\tilde{\sigma}_{i-1}^{A,T}$,
and $\tilde{\sigma}_{\mathrm{int}}$ are compatible in $\tilde{G}^{A,T}$.

\medskip We are now finally in a position to define a normal $5$-edge-coloring
$\tilde{\sigma}$ of $\tilde{G},$ and it is defined as follows:
\[
\left.  \tilde{\sigma}\right\vert _{M_{\mathrm{int}}}=\tilde{\sigma
}_{\mathrm{int}}=\left.  \sigma\right\vert _{M_{\mathrm{int}}}\text{ and
}\left.  \tilde{\sigma}\right\vert _{\mathcal{B}_{i}}=\tilde{\sigma}_{i}
\]
for $i=0,\ldots,g-1,$ where $\tilde{\sigma}_{i}$ is the coloring of
$\mathcal{B}_{i}$ (and the edge $u_{i}^{\prime}u_{i}^{\prime\prime}$ where
needed) as in Subcases 2.a, 2.b, and 2.c. In each of the three subcases, we
established that $\tilde{\sigma}_{i}$, $\tilde{\sigma}_{i-1}$, and
$\tilde{\sigma}_{\mathrm{int}}$ are compatible in $\tilde{G}$, so
$\tilde{\sigma}$ is a well-defined normal $5$-edge-coloring of $\tilde{G}$.
\end{proof}

Theorems \ref{Tm_oneSE} and \ref{Tm_twoSE} yield the following corollary.

\begin{corollary}
\label{Cor_evenSubgraph} Let $G$ be a snark with a normal $5$-edge-coloring
$\sigma,$ $C$ an even subgraph of $G,$ and $\tilde{G}\in\mathcal{G}
_{C}(\mathcal{A},\mathcal{B})$ a superposition of $G$. The superposition
$\tilde{G}$ has a normal $5$-edge-coloring if one of the following conditions holds:

\begin{itemize}
\item[i)] $\mathcal{B}_{i}$ is dock-right for every $e_{i}\in E(C)$;

\item[ii)] $\mathcal{B}_{i}$ is both doubly-right and doubly-left for every
$e_{i}\in E(C)$.
\end{itemize}
\end{corollary}

\begin{proof}
Since $C$ is an even subgraph of a cubic graph, it follows that $C$ is a
collection of vertex disjoint cycles in $G$. Notice that a normal
$5$-edge-coloring $\tilde{\sigma}$ constructed in the proofs of Theorems
\ref{Tm_oneSE} and \ref{Tm_twoSE} preserves the colors of $\sigma$ outside of
the considered cycle. Thus, a coloring $\sigma$ of $G$ can be extended to a
superposition along any vertex disjoint cycle independently.
\end{proof}

\begin{figure}[h]
\begin{center}
\includegraphics[scale=0.55]{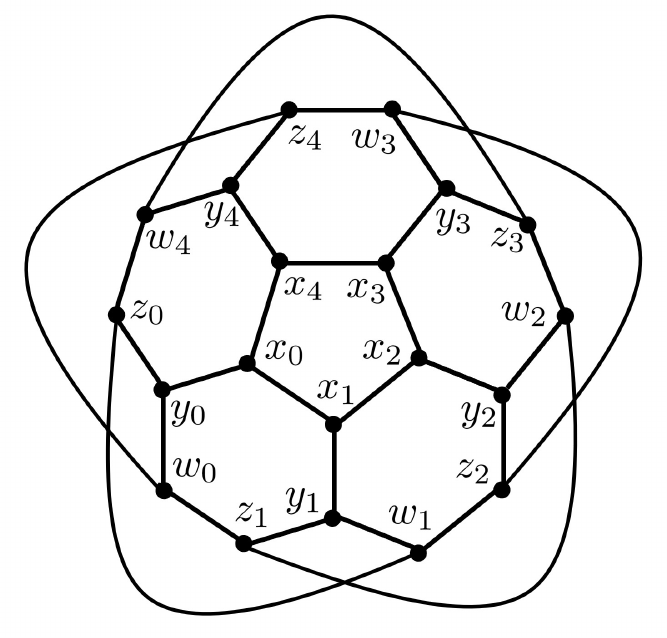}
\end{center}
\caption{The Flower snark $J_{5}.$}%
\label{Fig_flowerLabels}%
\end{figure}

\section{Application to snarks superpositioned by Flower snarks}

\begin{figure}[h]
\begin{center}%
\begin{tabular}
[t]{ll}%
a) & \raisebox{-0.9\height}{\includegraphics[scale=0.55]{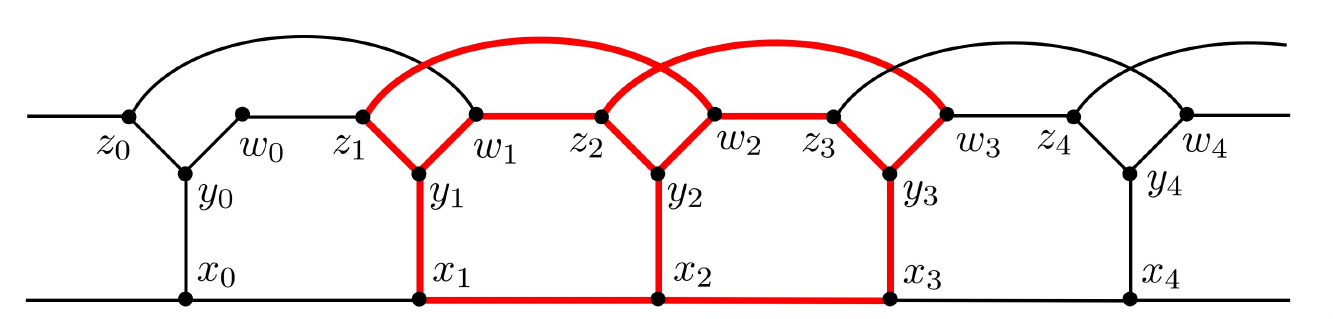}}\\
b) & \raisebox{-0.9\height}{\includegraphics[scale=0.55]{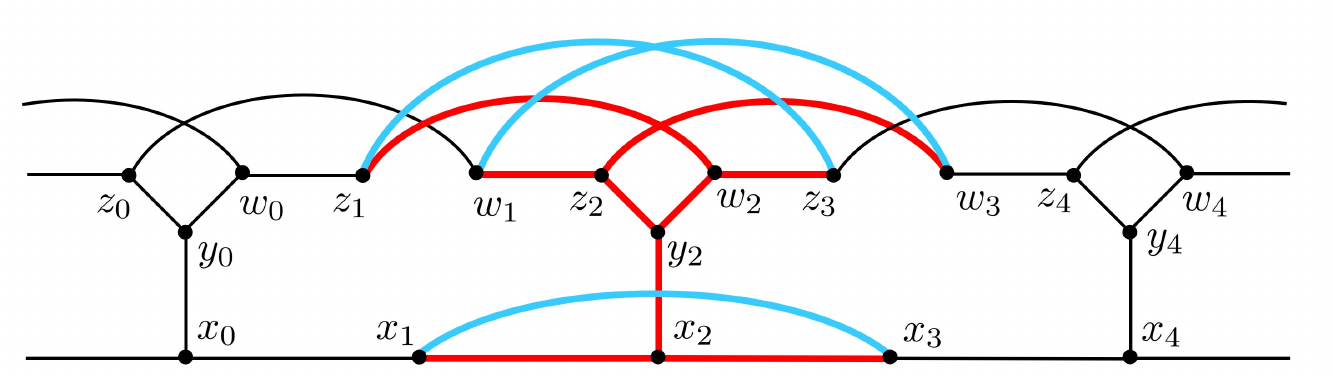}}\\
c) & \raisebox{-0.9\height}{\includegraphics[scale=0.55]{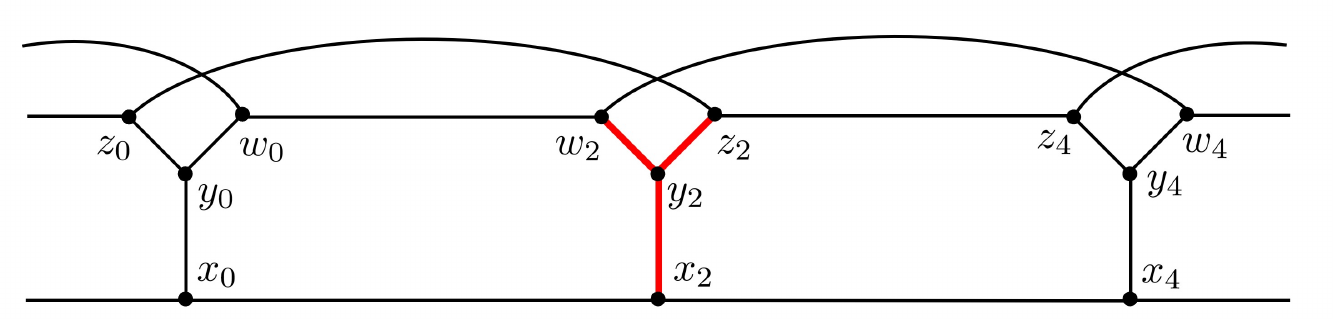}}
\end{tabular}
\end{center}
\caption{The figure shows: a) the graph $K$ in red as a subgraph of a larger
graph $G,$ b) $G$ with vertices $y_{1}$ and $y_{3}$ removed, and edges
$x_{1}x_{3},$ $z_{1}z_{3}$, $w_{1}w_{3}$ added which are the first two steps
of $K$-reduction of $G,$ c) $K$-reduced graph $G^{\prime}$ after the
contraction of triangles $x_{1}x_{2}x_{3},$ $z_{1}w_{2}z_{3}$, $w_{1}%
z_{2}w_{3}$ which is the third and last step of $K$-reduction.}%
\label{Fig_SteffenReduction}%
\end{figure}

We will now apply sufficient conditions of Theorems \ref{Tm_oneSE} and
\ref{Tm_twoSE} to snarks superpositioned by Flower snarks. Let us first define
the Flower snarks. A Flower snark $J_{r},$ for $r\geq5$ odd, is a graph with
the set of vertices $V(J_{r})=\cup_{i=0}^{r-1}\{x_{i},y_{i},z_{i},w_{i}\}$ and
the set of edges $E(J_{r})=\cup_{i=0}^{r-1}\{x_{i}x_{i+1},x_{i}y_{i},
y_{i}z_{i},y_{i}w_{i},z_{i}w_{i+1},w_{i}z_{i+1}\}$ where indices are taken
modulo $r$. The Flower snark for $r=5$ is illustrated by Figure
\ref{Fig_flowerLabels}. We will consider a superedge $H_{x,y}$ where $H=J_{r}$
and $x,y\in V(H)$ is any pair of vertices in $H$ at distance $\geq3.$ A
reduction introduced in \cite{Hagglund2014} makes it sufficient to consider
mostly $H=J_{5}$ and all pairs of vertices $x,y$ in $J_{5}$ at distance
$\geq3$ and only a few pairs from $J_{7}$. To be more precise, the graph $K$
shown in Figure \ref{Fig_SteffenReduction} in red and its reduction is
introduced by H\"{a}gglund and Steffen \cite{Hagglund2014}. If $K$ is a
subgraph of a larger cubic graph $G$, let $G^{\prime}$ be the graph obtained
from $G$ by:

\begin{itemize}
\item removing vertices $y_{1}$ and $y_{3},$

\item inserting edges $x_{1}x_{3},$ $z_{1}z_{3},$ and $w_{1}w_{3},$

\item contracting triangles $x_{1}x_{2}x_{3},$ $z_{1}w_{2}z_{3},$ and
$w_{1}z_{2}w_{3}.$
\end{itemize}

\noindent Such a graph $G^{\prime}$ is called a $K$\emph{-reduction} of $G$,
and $G$ is said to be $K$\emph{-reducible} to $G^{\prime}$. The following
theorem is established in \cite{Hagglund2014}.

\begin{theorem}
Let $G$ be a cubic graph that is $K$-reducible to a graph $G^{\prime}$. If
$G^{\prime}$ has a normal $5$-edge-coloring, then $G$ also has a normal $5$-edge-coloring.
\end{theorem}

In \cite{Hagglund2014}, it was also established that $J_{r},$ for $r\geq7,$ is
$K$-reducible to $J_{r-2}.$ This can be exploited in the case of a superedge
$(J_{r})_{x,y}$ in the following way. Due to the symmetry of $J_{r}$, the
vertex $x$ of the pair can be chosen from $\{x_{0},y_{0},z_{0},w_{0}\}$ and
the vertex $y$ from $\cup_{i=0}^{(r-1)/2}\{x_{i},y_{i},z_{i},w_{i}\}$, under
the condition that $x$ and $y$ are at distance $\geq3$. On the other hand, a
superedge $(J_{r})_{x,y},$ for $r\geq7,$ contains a subgraph $K$ induced by
vertices $\cup_{i=r-3}^{r-1}\{x_{i},y_{i},z_{i},w_{i}\}.$ More informally
said, the pair $x,y$ is always chosen in one half of the flower, and the other
half contains a subgraph $K.$ Consequently, any superposition which contains
such a superedge also contains a subgraph $K$ and it can be $K$-reduced. Thus,
we can talk of $K$-reduction of a superedge $(J_{r})_{x,y}$, for $r\geq7$, to
a superedge $(J_{r-2})_{x,y}.$ Moreover, the pair $x,y$ will remain at
distance $\geq3$ in $J_{r-2}$ in all cases, except for the case $r=7$ and only
a few pairs $(x,y)$ from it. To be more precise, we summarize all this in the
following formal observation.

\begin{remark}
Let $J_{r}$ be a Flower snark for $r\geq7,$ and $x,y$ a pair of vertices in
$J_{r}$ at distance $\geq3.$ A subgraph $K$ of $J_{r}$ can be chosen so that
the pair $x,y$ will remain at distance $\geq3$ in $K$-reduced $J_{r-2}$ in all
cases except in the case of $r=7$ and pairs $(x,y)\in\{(x_{0},x_{3}
),(z_{0},z_{3}),(z_{3},z_{0})\}$ up to a symmetry of $J_{7}.$
\end{remark}

Thus, to establish the existence of a normal $5$-edge-coloring for
superpositions by all Flower snarks it is sufficient to consider $J_{5}$ and
all pairs of vertices at the distance $\geq3$ in it, and additionally these
three pairs of vertices from $J_{7}.$ For the next proposition we need the
notion of a $2$-factor. So, a $2$\emph{-factor} of a subcubic graph $G$ is
defined as any spanning $2$-regular subgraph of $G$.

\begin{proposition}
\label{Prop_flowerDoubly}Let $J_{5}$ be a Flower snark and $x,y$ a pair of
vertices in $J_{5}$ at distance $\geq3.$ Then, a superedge $(J_{5})_{x,y}$ is
both doubly-right and doubly-left. The same holds for a superedge
$(J_{7})_{x,y}$ with $(x,y)\in\{(x_{0},x_{3}),(z_{0},z_{3}),(z_{3},z_{0})\}.$
\end{proposition}

\begin{proof}
Let us denote $H=J_{r}$ and $H^{\prime}=J_{r}-y.$ Denote by $x_{1},$ $x_{2}$
and $x_{3}$ the three neighbors of $x$ in $H.$ A $2$-factor $F$ of $H^{\prime
}$ is said to be \emph{right-good} if it contains precisely one odd-length
cycle $C_{x}$ for which $x\in V(C_{x})$ and all other cycles of $F$ are even.
For example, in each of the colorings from Figure \ref{Fig_completelyRight},
edges colored by $1$ and $2$ induce a right-good $2$-factor in $H.$

For each right-good $F,$ a corresponding right coloring of $H_{x,y}$ can be
obtained in the following way. Edges of $F$ are colored alternately by $1$ and
$2$ so that only $x$ is incident to a pair of edges of the same color, say
color $2.$ All edges of $H$ not belonging to $F$ are colored by $3.$ This
gives us a $3$-edge-coloring $\sigma$ of $H$ which is neither proper nor
normal. However, the restriction of $\sigma$ to $H_{x,y}$ is a right coloring
$R_{j}$, where $j$ denotes the left semiedge arising from the edge $x_{j}x$
not contained in $F.$ By $F_{j}$ we will denote a right-good $F$ such that
$C_{x}$ does not contain $x_{j}x.$ To establish that $H_{x,y}$ is doubly-left,
it is sufficient to show that there exists a right-good $F_{j}$ of $H_{x,y}$
for at least two distinct values of $j.$

Assuming vertex notation as in Figure \ref{Fig_flowerNumbers}, for $H=J_{5}$
this is up to symmetries done in Table \ref{tab_doublyRight5} in Appendix, and
for $H=J_{7}$ in Table \ref{tab_doublyRight7} also in Appendix. These two
tables consist of four columns. The first two give the pair of vertices $x$
and $y.$ The third column gives the vertex $x_{j},$ and the last column gives
the desired $2$-factor $F_{j}$ of the respective Flower snark.

To establish that $H_{x,y}$ is doubly-left, we first introduce the following
notions. Let $y_{1},$ $y_{2}$ and $y_{3}$ denote the neighbors of $y.$ A
$2$-factor $F$ of the snark $H$ is said to be \emph{left-good} if $F$ consists
of precisely two odd-length cycles $C_{x}$ and $C_{y}$ such that $x\in
V(C_{x})$ and $y\in V(C_{y}),$ and all other cycles of $F$ are even. For each
left-good $F,$ a corresponding left coloring of $H_{x,y}$ can be obtained in
the following way. Edges of $F$ are colored alternately by $1$ and $2$ so that
only $x$ and $y$ are incident to two edges of the same color, say $x$ of color
$1$ and $y$ of color $2$. All edges of $H$ not contained in $F$ are colored by
$3.$ This gives us a $3$-edge-coloring $\sigma$ of $H$ which is neither proper
nor normal. However, the restriction of $\sigma$ to $H_{x,y}$ is a left
coloring $L_{j,k}$, where $j$ denotes the left semiedge arising from the edge
$x_{j}x$ not contained in $F$ and $k$ denotes the right semiedge arising from
the edge $y_{k}y$ not contained in $F.$

By $F_{j,k}$ we will denote a left-good $F$ such that $C_{x}$ does not contain
$x_{j}x$ and $C_{y}$ does not contain $y_{k}y.$ To establish that $H_{x,y}$ is
doubly-left, it is sufficient to show that for every $j\in\{1,2,3\}$ there
exists a left-good $F_{j,k}$ for at least two distinct values of $k.$ Assuming
vertex notation as in Figure \ref{Fig_flowerNumbers}, for $H=J_{5}$ this is
done in Table \ref{tab_doublyLeft5} and for $H=J_{7}$ in Table
\ref{tab_doublyLeft7} in Appendix. These two tables each have five columns.
The first two columns give the pair $x$ and $y$ in the snark. The third and
the fourth columns give $x_{j}$ and $y_{k},$ respectively, and the last column
gives the desired $2$-factor $F_{j,k}$.
\end{proof}

The direct consequence of Proposition \ref{Prop_flowerDoubly} and Corollary
\ref{Cor_evenSubgraph} is the following result.

\begin{theorem}
\label{Tm_flower}Let $G$ be a snark which has a normal $5$-edge-coloring
$\sigma,$ $C$ an even subgraph of $G,$ and $\tilde{G}\in\mathcal{G}
_{C}(\mathcal{A},\mathcal{B})$ a superposition of $G$. If $\mathcal{B}_{i}
\in\{(J_{r})_{x,y}:x,y\in V(J_{r})$ and $d(x,y)\geq3\}$ for every $e_{i}\in
E(C)$, then $\tilde{G}$ has a normal $5$-edge-coloring.
\end{theorem}

Notice that according to Table \ref{tab_doublyRight5} from Appendix, a
superedge $(J_{5})_{x,y}$ is fully-right for
\[
(x,y)\in\{(y_{2},x_{0}),(y_{0},y_{1}),(y_{0},y_{2}),(y_{0},z_{2}),(y_{0}%
,w_{2})\}.
\]
So, to establish the existence of a normal $5$-edge-coloring for
superpositions with such superedges, even Theorem \ref{Tm_oneSE} suffices.
Theorem \ref{Tm_flower} implies the results of \cite{kineziEuropean}.

\section{Concluding remarks}

In this paper, we consider the Petersen Coloring Conjecture for a class of
superpositioned snarks. Since the Petersen Coloring Conjecture is equivalent
to stating that every bridgeless cubic graph has a normal $5$-edge-coloring,
we first introduce two sufficient conditions (in Theorems \ref{Tm_oneSE} and
\ref{Tm_twoSE}, respectively) for the existence of a normal $5$-edge-coloring
of superpositioned snarks $G_{C}(\mathcal{A},\mathcal{\mathcal{B}})$. These
snarks are obtained from any snark $G$ by choosing a cycle $C$ in $G$ and then
superpositioning each vertex of $C$ by a supervertex $\mathcal{A}_{i}
\in\{A,A^{\prime}\}$ and each edge of $C$ by a superedge $\mathcal{B}_{i}
\in\{H_{x,y}:H$ is a snark and $x,y\in V(H)$ such that $d(x,y)\geq3\}$. All
other vertices and edges of $G$ are superpositioned by themselves.

The sufficient condition of Theorem \ref{Tm_oneSE} is weaker than that of
Theorem \ref{Tm_twoSE}, but even so it applies to superpositions by infinitely
many distinct snarks. For example, it applies to all hypohamiltonian snarks
used as superedges but not to all possible ways of connecting them. In
particular, since Flower snarks are hypohamiltonian, Theorem \ref{Tm_oneSE}
implies that many snarks superpositioned by Flower snarks have a normal
$5$-edge-coloring. The condition of Theorem \ref{Tm_twoSE} is more demanding
and thus stronger. Its application to snarks superpositioned by Flower snarks
yields that all such superpositions have a normal $5$-edge-coloring.

It should be noted that for Flower snark superedges, the sufficient condition
of Theorem \ref{Tm_twoSE} is restated to involve the corresponding
$2$-factors, and then the existence of such $2$-factors is verified in silico.
An interesting question is whether the existence of the required $2$-factors
is due to the specific structure of Flower snarks, or the hypohamiltonian
property of such snarks. In the latter case, the results of this paper would
apply to all snarks superpositioned by hypohamiltonian snarks.

For a cubic graph $G,$ let $\mathrm{{NC}}(G)$ denote the set of all normal
$5$-colorings of $G$. The Petersen Coloring Conjecture is equivalent to the
claim that $\mathrm{{NC}}(G)\not =\emptyset$ for every bridgeless cubic graph
$G.$ Assuming that the Petersen Coloring Conjecture holds, i.e. that
$\mathrm{{NC}}(G)$ is indeed non-empty for every bridgeless cubic graph,
$\mathrm{{poor}}(G)$ is defined as the maximum number of poor edges among all
colorings from $\mathrm{{NC}}(G)$. Denote the Petersen graph by $P_{10},$ and
by $P_{10}^{\Delta}$ a graph obtained from $P_{10}$ by replacing one vertex with a triangle.
In \cite{SedSkrePaper1}, we have proposed the following conjecture.

\begin{conjecture}
\label{Conjecture_paper1}Let $G$ be a connected bridgeless cubic graph. If
$G\not =P_{10}$, then $\mathrm{{poor}}(G)>0$. Moreover, if $G\not =
P_{10},P_{10}^{\Delta}$, then $\mathrm{{poor}}(G)\geq6.$
\end{conjecture}

\noindent The results of this paper confirm this conjecture for the class of
snarks we considered. Namely, in a normal coloring constructed in the proof of
Theorem \ref{Tm_oneSE}, every superedge $\mathcal{B}_{i}$ is colored by a
(complementary) right coloring, and the right colorings of Flower snarks we
introduced use only three colors, so each edge in a right coloring of a
superedge is poor. Since a Flower snark superedge has at least 19 edges and
there are at least three superedges in a superposition, the claim of Conjecture
\ref{Conjecture_paper1} follows. For a normal coloring constructed in the
proof of Theorem \ref{Tm_twoSE}, at least every other superedge is colored by
a (complementary) right coloring, leading to the same conclusion. To be more
precise, in each of our constructions at least $\left\lfloor g/2\right\rfloor
\geq1$ superedges $\mathcal{B}_{i}$ of a superposition $\tilde{G}$ are colored
by $R_{j}^{i}$ or $\bar{R}_{j}^{i},$ and when $\mathcal{B}_{i}$ is a Flower
snark superedge then all its edges are poor in both $R_{j}^{i}$ and $\bar
{R}_{j}^{i}.$ Since every Flower snark superedge has at least $24$ edges, it
follows that every superposition by Flower snarks considered in this paper has
at least $24$ poor edges.

\bigskip

\bigskip\noindent\textbf{Acknowledgments.}~~Both authors acknowledge partial
support of the Slovenian research agency ARRS program\ P1-0383 and ARRS
project J1-3002. The first author also the support of Project
KK.01.1.1.02.0027, a project co-financed by the Croatian Government and the
European Union through the European Regional Development Fund - the
Competitiveness and Cohesion Operational Programme. The authors are also
thankful to the reviewers for the careful reading of the paper and valuable suggestions.

\newpage

\section{Appendix}

\begin{figure}[h]
\begin{center}%
\begin{tabular}
[c]{ll}%
\includegraphics[scale=0.55]{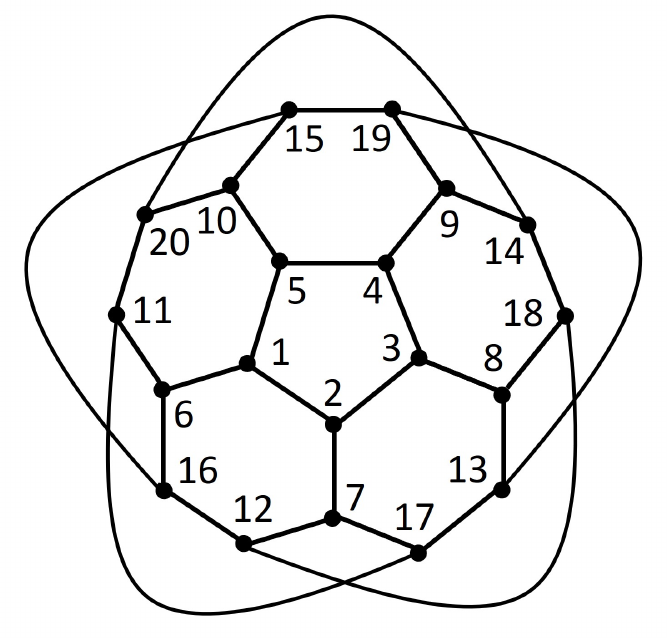} &
\includegraphics[scale=0.5]{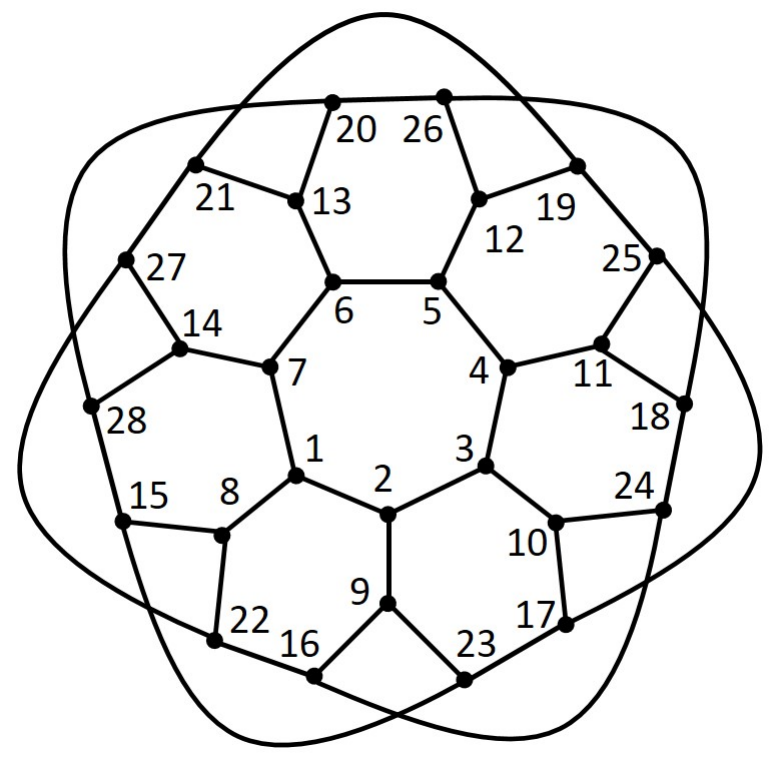}
\end{tabular}
\end{center}
\caption{The Flower snarks $J_{5}$ and $J_{7}$ with relabeled vertices.}%
\label{Fig_flowerNumbers}%
\end{figure}

We provide an illustration of relabeled snarks $J_{5}$ and $J_{7}$ in Figure
\ref{Fig_flowerNumbers}. This labeling is more convenient for giving all the
necessary $2$-factors of $J_{5}$ and $J_{7}$ which induce right and
left colorings mentioned in the proof of Proposition \ref{Prop_flowerDoubly}.
For this labeling of snarks $J_{5}$ and $J_{7},$ Table \ref{tab_doublyRight5}
has all $2$-factors of the graph $H-y,$ for $H=J_{5},$ which induce a
$j$-right coloring of a superedge $H_{x,y}$. Table \ref{tab_doublyRight7}
provides the same information only for $H=J_{7}.$ Tables \ref{tab_doublyLeft5}
and \ref{tab_doublyLeft7} contain $2$-factors of the graph $H=J_{5}$
and $H=J_{7},$ respectively, which induce a $j,k$-left coloring of a superedge
$H_{x,y}.$

%%%%%%%%%%%%%%%%%%%%%%%%%%%%%%%%%%%%%%%%%%%%%%%%%%%%%%%%%%%%%%%%%%%%%%%%%%%%%

\begin{center}
{\tiny \begin{longtable}{|l|l|l|l|}
		\caption{A $2$-factor of $H-y$ which induces a right coloring $R_j$ of $H_{x,y}$ for $H=J_5$.} \label{tab_doublyRight5} \\%
		
		\hline \multicolumn{1}{|c|}{\textbf{$x$}} & \multicolumn{1}{c|}{\textbf{$y$}} &
		\multicolumn{1}{c|}{\textbf{$x_j$}} & \multicolumn{1}{c|}{\text{$F_j$}} \\ \hline
		\endfirsthead
\text{{\tiny 1}} & \text{{\tiny 8}} & \text{{\tiny 2}} & \text{{\tiny (1, 5,
		4, 3, 2, 7, 17, 13, 19, 9, 14, 18, 12, 16, 15, 10, 20, 11, 6, 1)}}\\
\text{{\tiny 1}} & \text{{\tiny 8}} & \text{{\tiny 6}} & \text{{\tiny (1, 2,
		3, 4, 5, 1), (6, 11, 20, 10, 15, 16, 6), (7, 12, 18, 14, 9, 19, 13, 17, 7)}}\\
\text{{\tiny 8}} & \text{{\tiny 1}} & \text{{\tiny 3}} & \text{{\tiny (2, 3,
		4, 5, 10, 20, 14, 9, 19, 15, 16, 6, 11, 17, 13, 8, 18, 12, 7, 2)}}\\
\text{{\tiny 8}} & \text{{\tiny 1}} & \text{{\tiny 13}} & \text{{\tiny (2, 3,
		8, 18, 12, 16, 6, 11, 20, 14, 9, 4, 5, 10, 15, 19, 13, 17, 7, 2)}}\\
\text{{\tiny 8}} & \text{{\tiny 1}} & \text{{\tiny 18}} & \text{{\tiny (2, 3,
		8, 13, 17, 11, 6, 16, 15, 19, 9, 4, 5, 10, 20, 14, 18, 12, 7, 2)}}\\
\text{{\tiny 1}} & \text{{\tiny 12}} & \text{{\tiny 2}} & \text{{\tiny (8, 13,
		19, 9, 14, 18, 8), (1, 5, 4, 3, 2, 7, 17, 11, 20, 10, 15, 16, 6, 1)}}\\
\text{{\tiny 1}} & \text{{\tiny 12}} & \text{{\tiny 6}} & \text{{\tiny (3, 4,
		9, 14, 18, 8, 3), (1, 2, 7, 17, 13, 19, 15, 16, 6, 11, 20, 10, 5, 1)}}\\
\text{{\tiny 12}} & \text{{\tiny 1}} & \text{{\tiny 7}} & \text{{\tiny (2, 3,
		4, 5, 10, 15, 19, 9, 14, 20, 11, 6, 16, 12, 18, 8, 13, 17, 7, 2)}}\\
\hline
\end{longtable}}
\end{center}

%%%%%%%%%%%%%%%%%%%%%%%%%%%%%%%%%%%%%%%%%%%%%%%%%%%%%%%%%%%%%%%%%%%%%%%%%%%%%

\begin{center}
{\tiny \begin{longtable}{|l|l|l|l|}
			\caption{A $2$-factor of $H-y$ which induces a right coloring $R_j$ of $H_{x,y}$ for $H=J_7$.} \label{tab_doublyRight7} \\%
			
			\hline \multicolumn{1}{|c|}{\textbf{$x$}} & \multicolumn{1}{c|}{\textbf{$y$}} &
			\multicolumn{1}{c|}{\textbf{$x_j$}} & \multicolumn{1}{c|}{\text{$F_j$}} \\ \hline
			\endfirsthead
			
			\text{{\tiny 1}} & \text{{\tiny 4}} & \text{{\tiny 2}} & \text{{\tiny (2, 3,
					10, 17, 23, 9, 2), (1, 7, 14, 21, 22, 16, 24, 18, 11, 25, 19, 27, 13, 6, 5,
					12, 26, 20, 28, 15, 8, 1)}}\\
			\text{{\tiny 1}} & \text{{\tiny 4}} & \text{{\tiny 8}} & \text{{\tiny (5, 6,
					13, 27, 19, 12, 5), (1, 2, 3, 10, 24, 16, 9, 23, 17, 25, 11, 18, 26, 20, 28,
					15, 8, 22, 21, 14, 7, 1)}}\\
			\text{{\tiny 15}} & \text{{\tiny 18}} & \text{{\tiny 8}} & \text{{\tiny (12,
					19, 27, 13, 20, 26, 12), (2, 3, 10, 24, 16, 9, 2), (1, 7, 6, 5, 4, 11, 25, 17,
					23, 15, 28, 14, 21, 22, 8, 1)}}\\
			\text{{\tiny 15}} & \text{{\tiny 18}} & \text{{\tiny 23}} & \text{{\tiny (12,
					19, 27, 13, 20, 26, 12), (1, 7, 6, 5, 4, 11, 25, 17, 23, 9, 2, 3, 10, 24, 16,
					22, 21, 14, 28, 15, 8, 1)}}\\
			\text{{\tiny 18}} & \text{{\tiny 15}} & \text{{\tiny 11}} & \text{{\tiny (1,
					2, 3, 10, 17, 23, 9, 16, 24, 18, 26, 12, 19, 25, 11, 4, 5, 6, 7, 14, 28, 20,
					13, 27, 21, 22, 8, 1)}}\\
			\text{{\tiny 18}} & \text{{\tiny 15}} & \text{{\tiny 24}} & \text{{\tiny (9,
					16, 24, 10, 17, 23, 9), (1, 2, 3, 4, 5, 12, 19, 25, 11, 18, 26, 20, 28, 14, 7,
					6, 13, 27, 21, 22, 8, 1)}}\\ \hline
			
		\end{longtable}
}
\end{center}

%%%%%%%%%%%%%%%%%%%%%%%%%%%%%%%%%%%%%%%%%%%%%%%%%%%%%%%%%%%%%%%%%%%%%%%%%%%%%

\begin{center}
{\tiny \begin{longtable}{|l|l|l|l|l|}
			\caption{A $2$-factor of $H$ which induces a left coloring $L_{j,k}$ of $H_{x,y}$ for $H=J_5$.} \label{tab_doublyLeft5} \\%
			
			\hline \multicolumn{1}{|c|}{\textbf{$x$}} & \multicolumn{1}{c|}{\textbf{$y$}} &
			\multicolumn{1}{c|}{\textbf{$x_j$}} &
			\multicolumn{1}{c|}{\textbf{$y_k$}} &
			\multicolumn{1}{c|}{\text{$F_{j,k}$}} \\ \hline
			\endfirsthead
			
			\text{{\tiny 1}} & \text{{\tiny 8}} & \text{{\tiny 2}} & \text{{\tiny 3}} &
			\text{{\tiny (1, 5, 4, 3, 2, 7, 17, 11, 6, 1), (8, 13, 19, 9, 14, 20, 10, 15,
					16, 12, 18, 8)}}\\
			\text{{\tiny 1}} & \text{{\tiny 8}} & \text{{\tiny 2}} & \text{{\tiny 13}} &
			\text{{\tiny (1, 5, 4, 9, 19, 13, 17, 11, 6, 1), (2, 3, 8, 18, 14, 20, 10, 15,
					16, 12, 7, 2)}}\\
			\text{{\tiny 1}} & \text{{\tiny 8}} & \text{{\tiny 2}} & \text{{\tiny 18}} &
			\text{{\tiny (1, 5, 4, 9, 14, 18, 12, 16, 6, 1), (2, 3, 8, 13, 19, 15, 10, 20,
					11, 17, 7, 2)}}\\
			\text{{\tiny 1}} & \text{{\tiny 8}} & \text{{\tiny 5}} & \text{{\tiny 3}} &
			\text{{\tiny (1, 2, 3, 4, 5, 10, 20, 11, 6, 1), (7, 12, 16, 15, 19, 9, 14, 18,
					8, 13, 17, 7)}}\\
			\text{{\tiny 1}} & \text{{\tiny 8}} & \text{{\tiny 5}} & \text{{\tiny 13}} &
			\text{{\tiny (3, 4, 5, 10, 15, 16, 12, 18, 8, 3), (1, 2, 7, 17, 13, 19, 9, 14,
					20, 11, 6, 1)}}\\
			\text{{\tiny 1}} & \text{{\tiny 8}} & \text{{\tiny 5}} & \text{{\tiny 18}} &
			\text{{\tiny (3, 4, 5, 10, 20, 11, 17, 13, 8, 3), (1, 2, 7, 12, 18, 14, 9, 19,
					15, 16, 6, 1)}}\\
			\text{{\tiny 1}} & \text{{\tiny 8}} & \text{{\tiny 6}} & \text{{\tiny 3}} &
			\text{{\tiny (1, 2, 3, 4, 5, 1), (6, 11, 17, 7, 12, 18, 8, 13, 19, 9, 14, 20,
					10, 15, 16, 6)}}\\
			\text{{\tiny 1}} & \text{{\tiny 8}} & \text{{\tiny 6}} & \text{{\tiny 13}} &
			\text{{\tiny (1, 2, 7, 17, 13, 19, 15, 10, 5, 1), (3, 4, 9, 14, 20, 11, 6, 16,
					12, 18, 8, 3)}}\\
\hline
			
		\end{longtable}
}
\end{center}

%%%%%%%%%%%%%%%%%%%%%%%%%%%%%%%%%%%%%%%%%%%%%%%%%%%%%%%%%%%%%%%%%%%%%%%%%%%%%

\begin{center}
{\tiny \begin{longtable}{|l|l|l|l|l|}
			\caption{A $2$-factor of $H$ which induces a left coloring $L_{j,k}$ of $H_{x,y}$ for $H=J_7$.} \label{tab_doublyLeft7} \\%
			
			\hline \multicolumn{1}{|c|}{\textbf{$x$}} & \multicolumn{1}{c|}{\textbf{$y$}} &
			\multicolumn{1}{c|}{\textbf{$x_j$}} &
			\multicolumn{1}{c|}{\textbf{$y_k$}} &
			\multicolumn{1}{c|}{\text{$F_{j,k}$}} \\ \hline
			\endfirsthead
			
			\text{{\tiny 1}} & \text{{\tiny 4}} & \text{{\tiny 2}} & \text{{\tiny 3}} &
			\text{{\tiny (1, 7, 14, 28, 20, 26, 12, 19, 25, 17, 23, 15, 8, 1), (2, 3, 10,
					24, 18, 11, 4, 5, 6, 13, 27, 21, 22, 16, 9, 2)}}\\
			\text{{\tiny 1}} & \text{{\tiny 4}} & \text{{\tiny 2}} & \text{{\tiny 5}} &
			\text{{\tiny (1, 7, 14, 28, 20, 26, 18, 24, 10, 17, 23, 15, 8, 1), (2, 3, 4,
					11, 25, 19, 12, 5, 6, 13, 27, 21, 22, 16, 9, 2)}}\\
			\text{{\tiny 1}} & \text{{\tiny 4}} & \text{{\tiny 2}} & \text{{\tiny 11}} &
			\text{{\tiny (2, 3, 4, 5, 6, 13, 27, 21, 22, 16, 9, 2), (1, 7, 14, 28, 20, 26,
					12, 19, 25, 11, 18, 24, 10, 17, 23, 15, 8, 1)}}\\
			\text{{\tiny 1}} & \text{{\tiny 4}} & \text{{\tiny 7}} & \text{{\tiny 3}} &
			\text{{\tiny (4, 5, 6, 7, 14, 21, 22, 16, 9, 23, 17, 25, 11, 4), (1, 2, 3, 10,
					24, 18, 26, 12, 19, 27, 13, 20, 28, 15, 8, 1)}}\\
			\text{{\tiny 1}} & \text{{\tiny 4}} & \text{{\tiny 7}} & \text{{\tiny 5}} &
			\text{{\tiny (3, 4, 11, 18, 26, 20, 13, 27, 21, 22, 16, 24, 10, 3), (1, 2, 9,
					23, 17, 25, 19, 12, 5, 6, 7, 14, 28, 15, 8, 1)}}\\
			\text{{\tiny 1}} & \text{{\tiny 4}} & \text{{\tiny 7}} & \text{{\tiny 11}} &
			\text{{\tiny (3, 4, 5, 6, 7, 14, 21, 22, 16, 24, 10, 3), (1, 2, 9, 23, 17, 25,
					11, 18, 26, 12, 19, 27, 13, 20, 28, 15, 8, 1)}}\\
			\text{{\tiny 1}} & \text{{\tiny 4}} & \text{{\tiny 8}} & \text{{\tiny 3}} &
			\text{{\tiny (1, 2, 3, 10, 17, 25, 19, 12, 26, 20, 28, 14, 7, 1), (4, 5, 6,
					13, 27, 21, 22, 8, 15, 23, 9, 16, 24, 18, 11, 4)}}\\
			\text{{\tiny 1}} & \text{{\tiny 4}} & \text{{\tiny 8}} & \text{{\tiny 5}} &
			\text{{\tiny (1, 2, 9, 16, 24, 18, 26, 20, 28, 14, 7, 1), (3, 4, 11, 25, 19,
					12, 5, 6, 13, 27, 21, 22, 8, 15, 23, 17, 10, 3)}}\\
\hline
			
		\end{longtable}
}
\end{center}


\begin{thebibliography}{99}                                                                                               %


\bibitem {Adelson1973}G. M. Adelson-Velskii, V. K. Titov, On edge 4-chromatic
cubic graphs, \emph{Vopr. Kibern.} \textbf{1} (1973) 5--14.

\bibitem {BilkovaHana}H. B\'{\i}lkov\'{a}, Petersenovsk\'{e} obarven\'{\i} a
jeho varianty, Bachelor thesis, Charles University in Prague, Prague, 2012,
(in Czech).

\bibitem {Descartes1946}B. Descartes, Network-colourings, \emph{Math. Gaz.}
\textbf{32} (1948) 67--69.

\bibitem {MazzuccoloLupekhine}L. Ferrarini, G. Mazzuoccolo, V. Mkrtchyan,
Normal 5-edge-colorings of a family of Loupekhine snarks, \emph{AKCE Int. J.
Graphs Comb.} \textbf{17(3)} (2020) 720--724.

\bibitem {MazzuoccoloStephen}M. A. Fiol, G. Mazzuoccolo, E. Steffen, Measures
of edge-uncolorability of cubic graphs, \emph{Electron. J. Comb.}
\textbf{25(4)} (2018) P4.54.

\bibitem {Fiorini1983}S. Fiorini, Hypohamiltonian snarks, in: Graphs and Other
Combinatorial Topics, Proc. 3rd Czechoslovak Symp. on Graph Theory, Prague,
Aug. 24-27, 1982 (M. Fiedler, Ed.), Teubner-Texte zur Math., Bd. 59, Teubner,
Leipzig, 1983, 70--75.

\bibitem {Hagglund2014}J. H\"{a}gglund, E. Steffen, Petersen-colorings and
some families of snarks, \emph{Ars Math. Contemp.} \textbf{7} (2014) 161--173.

\bibitem {Isaacs1975}R. Isaacs, Infinite families of nontrivial trivalent
graphs which are not Tait colorable, \emph{Amer. Math. Monthly} \textbf{82}
(1975) 221--239.

\bibitem {Jaeger1988}F. Jaeger, Nowhere-zero flow problems, Selected topics in
graph theory, 3, Academic Press, San Diego, CA, 1988, 71--95.

\bibitem {Jaeger1985}F. Jaeger, On five-edge-colorings of cubic graphs and
nowhere-zero flow problems, \emph{Ars Comb.} \textbf{20-B} (1985) 229--244.

\bibitem {KocholSuperposition}M. Kochol, Snarks without small cycles, \emph{J.
Combin. Theory Ser. B} \textbf{67} (1996) 34--47.

\bibitem {Kochol2}M. Kochol, Superposition and constructions of graphs without
nowhere-zero k-flows, \emph{Eur. J. Comb.} \textbf{23} (2002) 281--306.

\bibitem {KocholFlower1}M. Kochol, A cyclically 6-edge-connected snark of
order 118, \emph{Discrete Math.} \textbf{161} (1996) 297--300.

\bibitem {kineziEuropean}S. Liu, R.-X. Hao, C.-Q. Zhang, Berge--Fulkerson
coloring for some families of superposition snarks, \emph{Eur. J. Comb.}
\textbf{96} (2021) 103344.

\bibitem {MacajovaSkovieraHypohamiltonian}E. M\'{a}\v{c}ajov\'{a}, M.
\v{S}koviera, Constructing Hypohamiltonian Snarks with Cyclic Connectivity 5
and 6, \emph{Electron. J. Comb.} \textbf{14} (2007), \#R18.

\bibitem {MacajovaRevisited}E. M\'{a}\v{c}ajov\'{a}, M. \v{S}koviera,
Superposition of snarks revisited, \emph{Eur. J. Comb.} \textbf{91} (2021) 103220.

\bibitem {Mazzucuolo2020normal}G. Mazzuoccolo, V. Mkrtchyan, Normal
edge-colorings of cubic graphs, \emph{J. Graph Theory} \textbf{94(1)} (2020) 75--91.

\bibitem {Mazzucuolo2020normal6}G. Mazzuoccolo, V. Mkrtchyan, Normal
6-edge-colorings of some bridgeless cubic graphs, \emph{Discret. Appl. Math.}
\textbf{277} (2020) 252--262.

\bibitem {Mkrtchjan2013}V. Mkrtchyan, A remark on the Petersen coloring
conjecture of Jaeger, \emph{Australas. J. Comb.} \textbf{56} (2013) 145--151.

\bibitem {NedelaSkovieraSurvey}R. Nedela, M. \v{S}koviera, Decompositions and
reductions of snarks, \emph{J. Graph Theory} \textbf{22} (1996) 253--279.

\bibitem {Riste2020}F. Pirot, J. S. Sereni, R. \v{S}krekovski, Variations on
the Petersen colouring conjecture, \emph{Electron. J. Comb.} \textbf{27(1)}
(2020) \#P1.8

\bibitem {Samal2011}R. \v{S}\'{a}mal, New approach to Petersen coloring,
\emph{Electr. Notes Discrete Math.} \textbf{38} (2011) 755--760.

\bibitem {SedSkrePaper1}J. Sedlar, R. \v{S}krekovski, Normal $5$-edge-coloring
of some snarks superpositioned by the Petersen graph, arXiv:2305.05981 [math.CO].
\end{thebibliography}
\end{document}